\documentclass[smallextended,referee,envcountsect,]{svjour3}
\smartqed
\usepackage{graphicx}
\usepackage{url}
\usepackage{cite}

\usepackage{enumerate}
\usepackage{enumitem}

\usepackage{bbm}
\usepackage{bm}  %%% command \bm{math expression} to bold the 'math expression'

\usepackage{amssymb}
\usepackage{amsmath}
\usepackage{cleveref}  %%% after \usepackage{amsmath}

\newcommand{\R}{\mathbb{R}}
\newcommand{\N}{\mathbb{N}}

\newcommand{\Linop}{\mathbb{L}}

\newcommand{\norm}[1]{\left\|#1\right\|}

\newcommand{\sign}{\mathrm{sign}}

\newcommand{\supp}{\mathrm{supp}}

\newcommand{\diam}{\mathrm{diam}}

\usepackage{todonotes}

\begin{document}

\title{Optimal control of a two-dimensional elliptic equation with exponential nonlinearity and Dirac measure data}

%\subtitle{Using  the  LaTex Template}
\author{Vu Huu Nhu}

%\author{Vu Huu Nhu\orcidlink{0000-0003-4279-3937}}

\institute{Vu Huu Nhu,  Corresponding author \at
			Faculty of Fundamental Sciences\\
			 PHENIKAA University\\
			 Yen Nghia, Ha Dong, Hanoi 12116, Vietnam\\
              nhu.vuhuu@phenikaa-uni.edu.vn
              %\\
              %\url{https://orcid.org/0000-0003-4279-3937}
%           \and
%              Virtual Author,  Corresponding author  \at
%              University of Dreamland \\
%              Dreamland\\
%              author2@example.com
}

\date{Received: date / Accepted: date}
%The correct dates will be entered by the editor.

\maketitle

\begin{abstract}
	This work addresses an optimal control problem for a semilinear elliptic equation in two-dimensional space, characterized by an exponential nonlinearity and a singular source term. The source is modeled as a finite linear combination of Dirac measures concentrated at a fixed set of distinct points. The control variable is a finite-dimensional vector whose components represent the masses assigned to these point sources.	
	Due to the interplay between the exponential nonlinearity and the singular measure data, the state equation is generally ill-posed and admits a unique very weak solution only when the largest component of the control vector does not surpass a certain critical threshold. Consequently, the control-to-state operator might be continuously differentiable only on an open subset of the control space.	
	To derive first-order optimality conditions for the original problem, we introduce a family of regularized problems by imposing box constraints on the control variables. These constraints are chosen such that the admissible control sets of the regularized problems lie entirely within the open subset where the control-to-state operator is smooth. By analyzing the optimality systems associated with the regularized problems and  passing to the limit, we obtain necessary optimality conditions for the original, unregularized problem.
\end{abstract}

%All acknowledgements should be placed in the back of the paper after Conclusions..

\section{Introduction}
\label{sec:introduction}

%%% Liouville type equation
The study of elliptic partial differential equations with exponential nonlinearities and singular data %, especially in two dimensions, 
is a rich area at the intersection of analysis, geometry, and mathematical physics. One important example is the Liouville-type equation that arises naturally in the context of prescribing Gaussian curvature, conformal metrics, and mean field models; see, e.g. \cite{ChenLin2015,BandleFlucher1996,Ni1982}.
In \cite{ChenLin2015}, Chen and Lin examined a generalized mean field equation on a compact Riemann surface, incorporating exponential nonlinearities and singular source terms represented by the Dirac delta function. 
%Their work is particularly notable for its topological degree analysis of solutions in the presence of singularities. 
%Meanwhile, in \cite{BandleFlucher1996}, Bandle and Flucher studied the geometric and variational structure of solutions to Liouville-type equations through the framework of conformal invariants. 

%%% Optimal control relating to measures & Inverse problems 
The formulation of the control problem in a measure space is motivated by the fact that the resulting optimal controls naturally exhibit sparsity, in the sense that they are supported on small subsets. This feature is especially advantageous in applications such as the optimal placement of sensors or actuators; see, e.g., \cite{CasasClasonKunisch2012,ClasonKunisch2012}, as well as in the context of inverse problems involving the identification of Dirac impulses; see, e.g., \cite{Meignen2019}. Although similar sparsity-promoting effects can be achieved by incorporating $L^1$-norm control costs, such formulations generally do not guarantee the existence of minimizers due to the lack of compactness in $L^1$-spaces. This limitation is effectively addressed by adopting a measure-space framework, where the control cost is defined in terms of measures.

For optimal control problems governed by elliptic PDEs within the measure-space framework, we refer to \cite{CasasClasonKunisch2012,CasasKunisch2014,ClasonKunisch2012,Hoppe2023}. In \cite{CasasKunisch2014}, Casas and Kunisch derived second-order necessary and sufficient optimality conditions, as well as results on the stability of solutions, for problems involving semilinear equations where the nonlinearity satisfies a polynomial growth condition. The work \cite{CasasClasonKunisch2012} by Casas et al. investigates first-order optimality conditions and the approximation of control problems governed by linear state equations, showing that optimal control measures can be effectively approximated by linear combinations of Dirac measures.

More recently, in \cite{Otarola2024}, a semilinear elliptic optimal control problem was considered, where the nonlinearity in the state equation is similar to that investigated in \cite{CasasKunisch2014}, and the forcing term is modeled as a linear combination of Dirac measures. In this work, the author established both first-order and second-order (necessary and sufficient) optimality conditions, along with error estimates for discrete approximations obtained via finite element methods.

%%%% Main problem 
%\medskip 
In this paper, we study an optimal control problem governed by a two-dimensional elliptic partial differential equation with an exponential nonlinearity. The source term in the state equation consists of the sum of an $L^p$-function with $p>1$, and a measure represented as a linear combination of  Dirac measures concentrated at fixed, distinct points. 
The presence of the exponential nonlinearity combined with the singular measure data makes the analysis of the state equation particularly challenging. In general, the existence of a weak solution to the state equation requires certain smallness conditions on the positive masses of the control, as large values can lead to blow-up or loss of well-posedness. 
%This balance between the strength of the point sources and the exponential growth of the nonlinearity is critical in ensuring the existence, uniqueness, and stability of solutions.
As a result, the control-to-state mapping is defined only on a subset of the control space  and is  continuously differentiable  only on  another open subset thereof; see \Cref{thm:F-diff-S} below.

%%% The aim of the paper
The aim of this paper is to derive optimality conditions for the considered optimal control problem. To this end, we adopt a regularization approach by introducing a family of regularized problems in which box constraints on the controls are additionally imposed. The upper bounds of these constraints are chosen to lie within the open subset on which the control-to-state mapping is continuously differentiable. 
This allows us to obtain the corresponding optimality systems for the regularized problems by using standard arguments. 
By subsequently passing to the limit in these systems, we obtain the desired optimality conditions for the original problem.
Our approach is similar in spirit to that used in \cite{Barbu1984,Tiba1990,MeyerSusu2017,Constantin2017,NeittaanmakiTiba1994}, where non-smooth state equations are approximated by smooth approximations.

% The outline is not required, but we show an example here.
The paper is organized as follows. In \Cref{sec:probem-setting-mainresults}, we present the problem setting and state the main results. \Cref{sec:elliptic-equation-measure} is devoted to auxiliary estimates for solutions to linear elliptic equations with Dirac measure data, as well as to the Gâteaux differentiability of the Nemytskii operator associated with the exponential nonlinearity. These results are applied in \Cref{sec:state-equation} to establish the continuous differentiability of the control-to-state mapping on an open subset of the control space. The proofs of the main results are provided in \Cref{sec:proof-mainresult}. The paper concludes with an appendix containing estimates for Green's function on the two-dimensional unit ball and a proof of the openness of certain subsets in the measure space and in a Sobolev space.

\medskip 

%%% notation
\paragraph*{Notation. } \,
The symbol $|v|$ means either the absolute value of $v$ if $v$ is a real number or the usual Euclidean norm when $v$ is a two-dimensional vector in $\R^2$. 
For $x$ belonging a Banach space $X$ and $\rho >0$, the symbol $B_{X}(x,\rho)$ stands for the open ball in $X$  of radius $\rho >0$ centered at $x$.
For a bounded domain $\Omega \subset \R^2$ with Lipschitz boundary, the symbol $\mathcal{M}(\Omega)$ stands for the space of all regular Borel measures in $\Omega$.
For ${\bm\omega} = (\omega_1, \omega_2, \ldots, \omega_k) \in \R^k$, by $\| \bm \omega \|_{1}$ and  $\|\bm \omega \|_{\infty}$, we denote, respectively, the $l^1$- and $l^\infty$-norms of $\bm\omega$;
by $\bm\omega^+$ and $\bm\omega^{-}$, we denote the positive and negative parts of $\bm\omega$, respectively, that is, $\bm{\omega}^{+} := (\omega_1^+, \omega_2^+, \ldots, \omega_k^+)$ and  $\bm{\omega}^{-} := (\omega_1^{-}, \omega_2^{-}, \ldots, \omega_k^{-})$, where $\omega_i^{+}$ and $\omega_i^{-}$ are  the positive and negative parts of $\omega_i$  and defined, respectively, by 
\[
\omega^{+}_{i} := \max\{ \omega_i, 0\} \quad \text{and} \quad  \omega^{-}_{i} := \max\{ -\omega_i, 0\} \quad \text{for } 1 \leq i \leq k.
\]
We also write
$
\bm\omega_{\max} := \max\{ \omega_i \mid 1 \leq i \leq k \}. 
$
In $\R^k$, with a slight abuse of notation, the symbol $\leq$ denotes the coordinatewise order, that is, for two $k$-dimensional vectors $\bm{\omega}$ and $\bm{\eta}$, the notation $\bm{\omega} \leq \bm{\eta}$ means that $\omega_i \leq \eta_i$ for all $1 \leq i \leq k$.
The symbol for a null vector in $\R^k$ is $\bm{0}$, that is, $\bm{0} = {(0,0,\ldots,0)} \in \R^k$,
%$
%	\bm{0} = \underbrace{(0,0,\ldots,0)}_{k},
%$
and the vector whose components are all equal to $1$ is denoted by $\bm{1}$. Similarly, we set $\bm{4\pi} := 4\pi \bm{1}$. 
%means that it is a vector in $\R^k$ and all of its component are equal to $4\pi$.
For a set $A \subset \R^n$, we denote by $\diam A$  its diameter, and $|A|$ its Lebesgue measure, provided $A$ is measurable.
The symbol $\sign(t)$ stands for the sign of the real number $t$, i.e., $\sign(t) =1$ if $t>0$, $\sign(t) = -1$ if $t <0$, and $\sign(t) = [-1,1]$ if $t=0$.
We use $C$ to denote a generic positive constant, whose value may vary between occurrences. When a dependence on a parameter is relevant, we also write, e.g., $C_t$ or $C(t)$ to indicate that the constant depends only on  the parameter $t$.
Finally, for two Banach spaces $X$ and $Y$, the notation $\Linop(X,Y)$ denotes the space of linear continuous operators from $X$ to $Y$.

%%% Optimal control problem with control constraints
\section{Problem setting and main results}
\label{sec:probem-setting-mainresults}
Let $\Omega \subset \R^2$ be a bounded domain with Lipschitz boundary $\Gamma$, 
%Let $\Omega \subset \R^2$ be a bounded domain with $C^1$-boundary $\Gamma$, 
and let $x_1$, $x_2$, $\ldots$, $x_k$ be fixed and distinct points in $\Omega$.
In this paper, we study optimal control problems in which the associated state equation is given by the following Dirichlet problem
\begin{equation}
	\label{eq:P-state}
	\left\{
	\begin{aligned}
		-\Delta y + (e^y - 1) & = f_0(x) + \sum_{i =1}^k \eta_i \delta_{x_i} && \text{in } \Omega, \\
		y &= 0 && \text{on } \Gamma
	\end{aligned}
	\right.
\end{equation}
for $f_0 \in L^p(\Omega)$ with $p>1$, $\eta_i \in \R$ and    $\delta_{x_i}$ denoting the Dirac measure concentrated at $x_i$,  $1 \leq i \leq k$.
We say that a function $y \in W^{1,q}_0(\Omega)$ for all $1 \leq q < 2$ is a (very) weak solution to \eqref{eq:P-state} if there holds
\begin{equation*}
	\left\{
	\begin{aligned}
		& e^y \in L^1(\Omega),\\
		& \int_\Omega  -\Delta \phi y + (e^y-1)\phi dx = \int_\Omega f_0 \phi dx + \sum_{i =1}^k \eta_i \phi(x_i) \quad \text{for all } \phi \in C^\infty_0(\Omega). 
	\end{aligned}
	\right.
\end{equation*}
By the results of Vazquez \cite{Vazquez1983} (see also~\cite[Thm.~4.17]{BresisMarcusPonce2007}, and \cite{Bartolucci2005} for the three-dimensional case), we obtain the following proposition concerning the existence and uniqueness of solutions to \eqref{eq:P-state}. Although the work of Vazquez \cite{Vazquez1983} considers equation \eqref{eq:P-state} on the entire space $\R^2$, the conclusions--and the structure of the proof--remain valid for bounded domains in two dimensions.
\begin{proposition}
	\label{prop:state-equation}
	The equation \eqref{eq:P-state} admits a weak solution $y$ if and only if $\bm{\eta} = (\eta_1, \eta_2,\ldots, \eta_k)$ satisfies
	\begin{equation}
		\label{eq:SE-existence-solution-assumption-0}
		\bm\eta_{\max} \leq 4\pi.
	\end{equation} 
	Moreover, \eqref{eq:P-state} has at most one solution.
\end{proposition}

We consider the optimal control problem 
\begin{equation}
	\label{eq:P-original}
	\tag{P}
	\min \left \{  J(y, \bm{\eta}) \, \middle| \, (y, \bm{\eta}) = (y, \eta_1, \eta_2, \ldots, \eta_k) \in L^2(\Omega) \times \R^k \, \text{satisfies \eqref{eq:P-state}} \right\},
\end{equation}
where objective functional is given by
\begin{equation}
	\label{eq:P-objective}
	J(y, \bm{\eta}) =  \frac{1}{2} \norm{y- y_d}^2_{L^2(\Omega)} + \kappa \sum_{i =1}^k |\eta_i|
\end{equation}
with $\kappa >0$ and $y_d \in L^2(\Omega)$.

We first establish the existence of minimizers for the problem \eqref{eq:P-original}.

\begin{theorem}
	\label{thm:minimizer-existence}
		There exists at least one global minimizer $(\bar y, \bar{ \bm{\eta}})$ of \eqref{eq:P-original}. Moreover, any minimizer $(\tilde y, \tilde{ \bm{\eta}})$ of \eqref{eq:P-original} satisfies
		\begin{equation}
			\label{eq:minimizer-notbigger-4pi}
			\bm{\tilde \eta}_{\max} \leq 4\pi.
		\end{equation}
\end{theorem}
{\it Proof.}   
	Let $y_0 \in H^1_0(\Omega)$ be a unique solution to \eqref{eq:P-state} corresponding to $\eta_i :=0$ for all $1 \leq i \leq k$.
	Since $ (y_0, \bm{0} ) \in L^2(\Omega) \times \R^k$ is an admissible point of \eqref{eq:P-original}, there holds
	\begin{equation} \label{eq:value-min-bounded}
		\alpha := \inf \left\{ J(y, \bm{\eta}) \,\middle| \, (y, \bm{\eta}) \, \text{satisfies \eqref{eq:P-state}}  \right\} \leq J(y_0, \bm{0}).
	\end{equation}
	Let $\{(y_n, \bm{\eta}_n )\}$ with $\bm{\eta}_n = (\eta_1^n, \eta_2^n, \ldots, \eta_k^n)$, $ n \geq 1$, be a minimizing sequence, that is, $(y_n, \bm{\eta}_n )$ satisfies \eqref{eq:P-state}	and
	\begin{equation}
		\label{eq:min-value}
		\alpha = \lim\limits_{n \to \infty} J(y_n,\bm{\eta}_n).
	\end{equation}
	Since $y_n$ is a weak solution to \eqref{eq:P-state}, we deduce from  \eqref{eq:SE-existence-solution-assumption-0} in \Cref{prop:state-equation} that
	\[
	\eta_i^n \leq 4 \pi \quad  \text{for all } 1 \leq i \leq k, n \geq 1.
	\]
	In view of \eqref{eq:value-min-bounded}, \eqref{eq:min-value}, and the definition of $J$, there exists an integer $n_0$ such that 
	\[
	\kappa \sum_{i=1}^k |\eta_i^n| \leq J(y_0,\bm{0}) + 1 \quad \text{for all} \quad n \geq n_0,
	\]
	which implies the existence of subsequence $\{{n_m} \}$ such that
	$
	\eta_i^{n_m} \to \eta_i^* 
	$ as $m \to \infty$
	for some $\eta_i^* \in \R$ for all $1 \leq i \leq k$. Obviously, there holds
	$
		\eta_i^* \leq   4 \pi 
	$ for all $1 \leq i \leq k$.
	Combining this with \Cref{prop:state-equation} shows that
	\eqref{eq:P-state} admits a unique solution 
	$y^*$ in $W^{1,q}_0(\Omega)$ for all $1 \leq q < 2$, associated with 
	 $\bm{\eta}^* := (\eta_1^*, \eta_2^*, \ldots, \eta_k^*)$.
	  
	We now prove that $(y^*, \bm{\eta}^*)$ is a global minimizer of \eqref{eq:P-original}. To this end, we employ the stability property of \eqref{eq:P-state}; see \eqref{eq:S-continuity-W1q} in \Cref{lem:control2state-oper-continuity} below. 
	Thanks to \eqref{eq:S-continuity-W1q}, one has
	\[
	y_{n_m} \to y^* \quad \text{in } W^{1,q}_0(\Omega) \quad \text{as } m  \to \infty
	\]
	for all $ 1 \leq q <2$. From this and the continuous embedding $W^{1,q}_0(\Omega) \hookrightarrow L^2(\Omega)$ for some $q$ sufficiently close to $2$, there holds
	$
	y_{n_m} \to y^* 
	$
	in $L^2(\Omega)$.
	There thus holds
	$
	\alpha = \lim\limits_{m \to \infty }J(y_{n_m}, \bm{\eta}_{n_m}) = J(y^*,\bm{\eta}^*). 
	$
	Consequently, $(y^*,\bm{\eta}^*)$ is a global solution to \eqref{eq:P-original}.
	Finally, as a direct consequence of \Cref{prop:state-equation},  \eqref{eq:minimizer-notbigger-4pi} is satisfied by any minimizer $(\tilde y, \tilde{ \bm{\eta}})$ of \eqref{eq:P-original}. 
\qed   

\begin{remark}
	\label{rem:restriction-feasible-set}
	As a result of \eqref{eq:SE-existence-solution-assumption-0} (and of \eqref{eq:minimizer-notbigger-4pi}), the feasible set of \eqref{eq:P-original} can be reduced to 
	\[
	\left\{(y, \bm{\eta}) \,\middle| \, \bm{\eta}_{\max} \leq 4 \pi \, \text{and $(y, \bm{\eta})$ satisfies \eqref{eq:P-state}}  \right\}.
	\]
\end{remark}

Since the upper bound $4\pi$ is critical (see assertion \ref{item:4pi} in \Cref{prop:critical-4pi}), we must carefully analyze the indices at which  components of a control minimizer attain the value $4\pi$. From now on, for any $\bm{\eta} \in \R^k$, we define the following critical index set 
\begin{equation}
	\label{eq:indices-4pi}
	I_{4\pi}(\bm{ \eta}) := \{ i \in \N \mid 1 \leq i \leq k, \eta_i = 4\pi \}.
\end{equation}

We now state the main result of this paper. %Its proof is left to \Cref{sec:proof-mainresult}.
\begin{theorem}
	\label{thm:1st-OCs-P-orig}
	Assume that $(\bar y, \bm{\bar \eta})$ is a local minimizer of \eqref{eq:P-original}. 
%	Further, assume that it is a locally strict minimizer of \eqref{eq:P-original}, provided that $I_{4\pi}(\bm{\bar \eta}) \neq \emptyset$.
	Assume further that one of the following conditions is fulfilled:
%	\begin{align}
%		& I_{4\pi}(\bm{\bar \eta}) = \emptyset; \, \text{or} \label{ass:regularcase}\\
%		& I_{4\pi}(\bm{\bar \eta}) \neq \emptyset \, \text{and} \, (\bar y, \bm{\bar \eta}) \, \text{is a locally strict minimizer of \eqref{eq:P-original}.} \label{ass:irrergularcase}
%	\end{align}
	\begin{enumerate}[label=(A\arabic*),leftmargin=1cm]
		\item \label{ass:regularcase}
		$I_{4\pi}(\bm{\bar \eta}) = \emptyset$; or
		
		\item \label{ass:irrergularcase}
		$I_{4\pi}(\bm{\bar \eta}) \neq \emptyset$ and $(\bar y, \bm{\bar \eta})$ is a strict local  minimizer of \eqref{eq:P-original}.
	\end{enumerate}	
	Then,  there exists a $\bar\varphi \in H^1_0(\Omega) \cap L^\infty(\Omega)$ that, together with $\bm{\bar \eta}$ and $\bar y$, satisfies the following  optimality conditions:
	\begin{subequations}
		\label{eq:OCs-1st-P-orig}
		\begin{align}
			& 
			-\Delta \bar y + (e^{\bar y} - 1)  = f_0(x) + \sum_{i =1}^k \bar\eta_i \delta_{x_i} \, \text{in } \Omega, \quad
			\bar y = 0 \, \text{on } \Gamma, \label{eq:OC-1st-state-P-orig}\\
%			
%			\left\{
%			\begin{aligned}
%				-\Delta \bar y + (e^{\bar y} - 1)  &= f_0(x) + \sum_{i =1}^k \bar\eta_i \delta_{x_i} && \text{in } \Omega, \\
%				\bar y &= 0 && \text{on } \Gamma,
%			\end{aligned}
%			\right.  \label{eq:OC-1st-state-P-orig} \\
			& 
			-\Delta \bar \varphi + e^{\bar y} \bar\varphi = \bar y - y_d \, \text{in } \Omega, \quad
			\bar \varphi = 0 \, \text{on } \Gamma, \label{eq:OC-1st-adjoint-state-P-orig} 
%			
%			\left\{
%			\begin{aligned}
%				-\Delta \bar \varphi + e^{\bar y} \bar\varphi &= \bar y - y_d && \text{in } \Omega, \\
%				\bar \varphi &= 0 && \text{on } \Gamma,
%			\end{aligned}
%			\right. \label{eq:OC-1st-adjoint-state-P-orig} 
			\intertext{and}
			& \left\{
			\begin{aligned}
				%& - \frac{\bar\varphi(x_i)}{\kappa} \geq 1 && \text{if } \bar\eta_i = M_i, \\
				& - \frac{\bar\varphi(x_i)}{\kappa} = 1 && \text{if } \bar\eta_i \in (0, 4\pi), \\ 
				& - \frac{\bar\varphi(x_i)}{\kappa} \in [-1,1] && \text{if } \bar\eta_i = 0, \\
				& - \frac{\bar\varphi(x_i)}{\kappa} =- 1 && \text{if } \bar\eta_i <0.
			\end{aligned}
			\right.
			\quad \text{for } 1 \leq i \leq k. 
			\label{eq:OC-1st-variational-P-orig} 
		\end{align}
	\end{subequations}
	Moreover, $\bar\varphi$ is continuous at $x_i$ for all $i \notin I_{4\pi}(\bm{\bar \eta})$.
\end{theorem}
The proof of this main result is deferred to \Cref{sec:proof-mainresult}, where assumption \ref{ass:regularcase} is considered in \Cref{subsec:proof-main-regular}, and hypothesis \ref{ass:irrergularcase} is investigated in \Cref{subsec:proof-main-irregular}.
\begin{remark}
	In \eqref{eq:OC-1st-variational-P-orig} of \Cref{thm:1st-OCs-P-orig}, we do not state any condition at the critical index $i$, where $\bar \eta_i = 4\pi$. This is due to the lack of $L^{1 + \tau}$-regularity of $e^{\bar y}$ in a vicinity of $x_i$ for some $\tau >0$; see, assertion \ref{item:4pi} in \Cref{prop:critical-4pi} below. 
\end{remark}
\section{Auxiliary results on elliptic equations involving measures and the exponential Nemytskii operator}
\label{sec:elliptic-equation-measure}

We begin this section by deriving estimates for solutions to Poisson's equation with a source being a measure concentrated at a set of finitely many points.
This result generalizes the one in \cite[Thm.~1]{BrezisMerle1991}.
\begin{proposition}
	\label{prop:Poisson-exponential-esti-many-Dirac}
	Assume that $\bm{\omega} = (\omega_1, \omega_2, \ldots, \omega_k) \in \R^k$ such that $\omega_i >0$ for all $i \in \{1,2,\ldots, k\}$. 
	%Let  $x_1, x_2, \ldots, x_k$ be mutually distinct points in $G \subset \R^2$. 
	Let $y \in W^{1,q}_0(\Omega)$, $1 \leq q <2$, be a unique solution to 
	\begin{equation}
		\label{eq:Poisson-many-Dirac}
		-\Delta y  = \sum_{i=1}^k \omega_i \delta_{x_i} \, \text{in } \Omega, \quad
		y = 0 \, \text{on } \Gamma,
%		
%		\left\{
%		\begin{aligned}
%			-\Delta y & = \sum_{i=1}^k \omega_i \delta_{x_i} && \text{in } \Omega, \\
%			y &= 0 && \text{on } \Gamma.
%		\end{aligned}
%		\right.
	\end{equation}
	and let $r_0 > 0$  be constant such that balls $B_{\R^2}(x_i, r_0) \subset \Omega$, $i=1,2,\ldots, k,$ are pairwise disjoint. 
	%For $k=1$, we take $r_0 = \frac{1}{2} \diam \Omega$.
	Then, for any $\alpha \in (0, 4 \pi)$, there hold
	\begin{multline}
		\label{eq:exponential-esti-Dirac-many}
		\int_\Omega \exp [\frac{(4\pi - \alpha)|y(x)|}{\bm\omega_{\max}} ] dx \leq \pi (\frac{2R}{r_0})^{(2 -\frac{\alpha}{2 \pi})\sum_{i=1}^k \frac{\omega_i}{\bm\omega_{\max} }} \\
		\begin{aligned}				
			 & \times [ (R^2 - kr_0^2) + 2r_0^2 \sum_{i=1}^{k} ( 2 - (2 - \frac{\alpha}{2 \pi}) \frac{\omega_i}{\bm\omega_{\max} } )^{-1} ] 
		\end{aligned}
	\end{multline}
	and 
	\begin{equation}
		\label{eq:exponential-esti-Dirac-many-Bxj}
		\int_{B_{\R^2}(x_j, r_0)} \exp [\frac{(4\pi - \alpha)|y(x)|}{\omega_j} ] dx \leq \frac{4\pi^2 r_0^2}{\alpha}  (\frac{2R}{r_0})^{(2 -\frac{\alpha}{2 \pi}) \sum_{i=1}^k \frac{\omega_i}{\omega_j }}
	\end{equation}
	for all $1 \leq j \leq k$ with $R := \frac{1}{2} \diam \Omega$.
	% Here $y$ can be extended by zero outside $\Omega$.
\end{proposition}
%\begin{remark}
%	\label{rem:exponential-esti}
%	When $k=1$ and thus $r_0 = R = \frac{1}{2} \diam \Omega$, the exponential estimate \eqref{eq:exponential-esti-Dirac-many} reduces to
%	\[
%	\int_\Omega \exp [\frac{(4\pi - \alpha)|y(x)|}{\omega_{1}} ] dx \leq \pi 2^{(2 -\frac{\alpha}{2 \pi})} \times 2R^2 \times ( \frac{\alpha}{2 \pi} )^{-1} = 	\frac{4 \pi^2}{\alpha}  (\diam \Omega)^2  2^{- \frac{\alpha}{2 \pi}}.
%	\] 
%	This extends the result in \cite[Thm.~1]{BrezisMerle1991}. 
%	%This shows \eqref{eq:exponential-esti-Dirac-single}.
%\end{remark}
%
%\medskip
%\noindent\textit{Proof of \Cref{prop:Poisson-exponential-esti-many-Dirac}}.
{\it Proof.}
By using a linear shift if necessary, we assume that $\Omega \subset B_{\R^2}(0,R) =: B_R$.
Let $\epsilon \in (0, \frac{r_0}{2})$ be arbitrary but fixed, and let $y_i^\epsilon$ and $\tilde{y}_i^\epsilon$, $1 \leq i \leq k$, be unique solutions to
	\begin{align}
		& 	\label{eq:Poisson-delta-i-G-BR}
		\left\{
		\begin{aligned}
			-\Delta y_i^\epsilon & = \phi_\epsilon(\cdot -x_i) && \text{in } \Omega, \\
			y_i^\epsilon &= 0 && \text{on } \Gamma
		\end{aligned} 
		\right.
		\quad \text{and}
		\quad
		\left\{
		\begin{aligned}
			-\Delta \tilde y_i^\epsilon & = \phi_\epsilon(\cdot -x_i) && \text{in } B_R, \\
			\tilde y_i^\epsilon &= 0 && \text{on } \partial B_R,
		\end{aligned}
		\right.
	\end{align}
respectively,
where $\{\phi_\epsilon\}_{0 < \epsilon < r_0/2}$ is a family of mollifiers satisfying
\[
0 \leq \phi_\epsilon \in C^\infty_c(\R^2), \quad \supp \phi_\epsilon \subset \overline{B_{\R^2}(0, \epsilon)}, \quad \int_{\R^2} \phi_\epsilon dx =1.
\]
A weak form of the maximum principle; see, e.g. \cite[Prop.~4.B.1]{BresisMarcusPonce2007}, therefore gives $y_i^\epsilon \geq 0$ a.e. in $\Omega$ and  $\tilde y_i^\epsilon \geq 0$ a.e. in $B_R$ for all $1 \leq i \leq k$. 
Again, using the maximum principle implies that
\begin{equation}
	\label{eq:yi-epsilon-comparison}
	0 \leq y_i^\epsilon \leq \tilde{y}_i^\epsilon \quad \text{a.e. in } \Omega.
\end{equation}
On the other hand, since $\phi_\epsilon(\cdot - x_i)$ converges weakly-star to $\delta_{x_i}$ in $\mathcal{M}(\Omega)$ as $\epsilon \to 0$, $y_i^\epsilon$ strongly converges to $y_i$ in $W^{1,q}_0(\Omega)$ that uniquely solves
\[
	-\Delta y_i  = \delta_{x_i} \, \text{in } \Omega, \quad
	y_i = 0 \, \text{on } \Gamma
%
%\left\{
%\begin{aligned}
%	-\Delta y_i & = \delta_{x_i} && \text{in } \Omega, \\
%	y_i &= 0 && \text{on } \Gamma
%\end{aligned} 
%\right.
\]
for all $1 \leq i \leq k$;
see, e.g. \cite[Thm.~2.1]{CasasKunisch2014}. 
Setting now ${y}^\epsilon := \sum_{i=1}^k \omega_i y_i^\epsilon$ and  $\tilde{y}^\epsilon := \sum_{i=1}^k \omega_i \tilde y_i^\epsilon$ thus yields 
\begin{equation}
	\label{eq:yi-epsilon-limit}
	\left\{
	\begin{aligned}
		& 	{y}^\epsilon \to y \quad \text{in } W^{1,q}_0(\Omega) \quad \text{and a.e. in } \Omega,\\
		&  \tilde{y}^\epsilon \geq {y}^\epsilon \geq 0 \quad \text{a.e. in } \Omega,
	\end{aligned}
	\right.
\end{equation}
where the estimates in \eqref{eq:yi-epsilon-limit} is due to \eqref{eq:yi-epsilon-comparison}.

\medskip 
Fixing $\alpha \in (0, 4 \pi)$ and $j \in \{1,2, \ldots, k\}$, we now estimate $\exp[\frac{(4 \pi - \alpha)\omega_i}{\bm\omega_{\max}}|\tilde{y}_i^\epsilon(x)|]$ and $\exp[\frac{(4 \pi - \alpha)\omega_i}{\omega_{j}}|\tilde{y}_i^\epsilon(x)|]$, $1 \leq i \leq k$. To this end, we first employ the Green representation formula; see, e.g. \cite[Cor.~2.2]{VeronVivier2000}, \cite[Thm.~12, Chap.~2]{Evans2010}, and \cite{ChipotQuittner2004}, for the second equation in \eqref{eq:Poisson-delta-i-G-BR}  to have
\begin{equation}
	\label{eq:Green-repre}
	\tilde{y}_i^\epsilon(x) = \int_{B_R} G(\frac{x}{R}, \frac{x'}{R})\phi_\epsilon(x' - x_i)dx' =  \int_{B_{\R^2}(x_i, \epsilon)} G(\frac{x}{R}, \frac{x'}{R})\phi_\epsilon(x' - x_i)dx' 
\end{equation}
{for a.e. } $x \in B_R$,
where $G$ denotes the Green function of the Laplace operator $-\Delta$ in the unit ball in $\R^2$. 
By \Cref{lem:Green-esti-unitball}, there holds
\begin{equation}
	\label{eq:yi-epsilon-esti-pointwise}
	0 \leq \tilde{y}_i^\epsilon(x) \leq \frac{1}{2 \pi} \int_{B_{\R^2}(x_i, \epsilon)} \ln ( \frac{2R}{|x-x'|} )  \phi_\epsilon(x' - x_i)dx' \quad \text{for a.e. } x \in B_R.
\end{equation}
We thus consider the following situations:

%\begin{itemize}
%\item 
\noindent$\bullet$
For a.e.  $x \in B_R \backslash \overline{ B_{\R^2}(x_i, r_0)}$, there holds $|x - x'| \geq r_0-\epsilon $ for all $x' \in B_{\R^2}(x_i, \epsilon)$. Then, one has
\[
0 \leq \tilde{y}_i^\epsilon(x) \leq \frac{1}{2 \pi} \int_{B_{\R^2}(x_i, \epsilon)} \ln (\frac{2R}{r_0 -\epsilon}) \phi_\epsilon(x' - x_i)dx' = \frac{1}{2 \pi}\ln \frac{2R}{r_0 -\epsilon}.
\]
This consequently shows that
\begin{equation}
	\label{eq:yi-epsilon-esti-pointwise-outside}
	\exp[\frac{(4 \pi - \alpha)\omega_i}{\bm\omega_{\max}}|\tilde{y}_i^\epsilon(x)|] \leq  ( \frac{2R}{r_0 -\epsilon} )^{(2 - \frac{\alpha}{2 \pi}) \frac{\omega_i}{\bm\omega_{\max}} }
\end{equation}
and
\begin{equation}
	\label{eq:yi-epsilon-esti-pointwise-outside-Bxj}
	\exp[\frac{(4 \pi - \alpha)\omega_i}{\omega_j}|\tilde{y}_i^\epsilon(x)|] \leq  ( \frac{2R}{r_0 -\epsilon} )^{(2 - \frac{\alpha}{2 \pi})  \frac{\omega_i}{\omega_j}}
\end{equation}
for a.e. $x \in B_R \backslash \overline{  B_{\R^2}(x_i, r_0)}$ and for all $1 \leq i \leq k$.

%%\item 
\noindent$\bullet$
For a.e. $x \in \overline{  B_{\R^2}(x_i, r_0)}$, in view of \eqref{eq:yi-epsilon-esti-pointwise} and the fact that $\supp \phi_\epsilon \subset \overline{B_{\R^2}(0, \epsilon)}$, there holds 
\begin{multline*}
	\exp[\frac{(4 \pi - \alpha)\omega_i}{\bm\omega_{\max}}|\tilde{y}_i^\epsilon(x)|] %\\
	\leq \exp [ (2 - \frac{\alpha}{2 \pi}) \frac{\omega_i}{\bm\omega_{\max}} \int_{B_{\R^2}(x_i, \epsilon)} \ln ( \frac{2R}{|x-x'|} )  \phi_\epsilon(x' - x_i)dx' ].
\end{multline*}
Applying now Jensen's inequality for integrals; see, .e.g. \cite[Thm.~B.50]{Leoni2017}, and using the identity 
$
\int_{B_{\R^2}(x_i, \epsilon)}\phi_\epsilon(x' - x_i)dx' =1,
$
we arrive at
\begin{multline*}
	\exp[\frac{(4 \pi - \alpha)\omega_i}{\bm\omega_{\max}}|\tilde{y}_i^\epsilon(x)|] \leq 
	\int_{B_{\R^2}(x_i, \epsilon)} \exp [ (2 - \frac{\alpha}{2 \pi}) \frac{\omega_i}{\bm\omega_{\max}} \ln ( \frac{2R}{|x-x'|} ) ]   \phi_\epsilon(x' - x_i) dx'\\
	\begin{aligned}
	%	& \leq 
	%	\int_{B_{\R^2}(x_i, \epsilon)} \exp [ (2 - \frac{\alpha}{2 \pi}) \frac{\omega_i}{\bm\omega_{\max}} \ln ( \frac{2R}{|x-x'|} ) ]   \phi_\epsilon(x' - x_i) dx' \\
		& = 
		\int_{B_{\R^2}(x_i, \epsilon)} ( \frac{2R}{|x-x'|} )^{(2 - \frac{\alpha}{2 \pi}) \frac{\omega_i}{\bm\omega_{\max}} }   \phi_\epsilon(x' - x_i) dx'.
	\end{aligned}
\end{multline*}
Integrating over ${ B_{\R^2}(x_i, r_0)}$ gives
\begin{multline*}
	\int_{ B_{\R^2}(x_i, r_0)} \exp[\frac{(4 \pi - \alpha)\omega_i}{\bm\omega_{\max}}|\tilde{y}_i^\epsilon(x)|]  dx \\
	\leq \int_{ B_{\R^2}(x_i, r_0)} dx \int_{B_{\R^2}(x_i, \epsilon)} ( \frac{2R}{|x-x'|} )^{(2 - \frac{\alpha}{2 \pi}) \frac{\omega_i}{\bm\omega_{\max}} }   \phi_\epsilon(x' - x_i) dx'.
\end{multline*}
A method of  integration by substitution then yields
\begin{multline*}
	\int_{ B_{\R^2}(x_i, r_0)} \exp[\frac{(4 \pi - \alpha)\omega_i}{\bm\omega_{\max}}|\tilde{y}_i^\epsilon(x)|]  dx  \\
	\begin{aligned}
		& \leq \int_{ B_{\R^2}(0,  r_0 + \epsilon)} d \eta \int_{B_{\R^2}(0, \epsilon)} ( \frac{2R}{|\eta|} )^{(2 - \frac{\alpha}{2 \pi}) \frac{\omega_i}{\bm\omega_{\max}} }   \phi_\epsilon(\psi) d\psi \\
		& =  \int_{ B_{\R^2}(0,  r_0 + \epsilon)} ( \frac{2R}{|\eta|} )^{(2 - \frac{\alpha}{2 \pi}) \frac{\omega_i}{\bm\omega_{\max}} } d \eta 
		 = 2 \pi \int_0^{r_0 + \epsilon }  ( \frac{2R}{r} )^{(2 - \frac{\alpha}{2 \pi}) \frac{\omega_i}{\bm\omega_{\max}} } rdr \\
		& = 2 \pi (2R)^{(2 - \frac{\alpha}{2 \pi}) \frac{\omega_i}{\bm\omega_{\max}} } \times  \frac{1}{2 - {(2 - \frac{\alpha}{2 \pi}) \frac{\omega_i}{\bm\omega_{\max}} }}(r_0 + \epsilon)^{2 -{(2 - \frac{\alpha}{2 \pi}) \frac{\omega_i}{\bm\omega_{\max}} }}.
	\end{aligned}
\end{multline*}
In conclusion, we have
\begin{multline}
	\label{eq:yi-epsilon-esti-inball}
	\int_{ B_{\R^2}(x_i, r_0)} \exp[\frac{(4 \pi - \alpha)\omega_i}{\bm\omega_{\max}}|\tilde{y}_i^\epsilon(x)|]  dx 
	\leq \frac{2 \pi (\epsilon + r_0)^2}{2 - {(2 - \frac{\alpha}{2 \pi}) \frac{\omega_i}{\bm\omega_{\max}} }}  ( \frac{2R}{\epsilon + r_0} )^{(2 - \frac{\alpha}{2 \pi}) \frac{\omega_i}{\bm\omega_{\max}} }
\end{multline}
for all $1 \leq i \leq k$.
Similarly, there holds for all $1 \leq i \leq k$ that
\begin{multline}
	\label{eq:yi-epsilon-esti-inball-Bxj}
	\int_{ B_{\R^2}(x_i, r_0)} \exp[\frac{(4 \pi - \alpha)\omega_i}{\omega_i}|\tilde{y}_i^\epsilon(x)|]  dx %\\
	\leq \frac{2 \pi (\epsilon + r_0)^2}{\frac{\alpha}{2 \pi} }  ( \frac{2R}{\epsilon + r_0} )^{(2 - \frac{\alpha}{2 \pi})  }.
\end{multline}
%\end{itemize}

\medskip 
We have shown the estimates \eqref{eq:yi-epsilon-esti-pointwise-outside}--\eqref{eq:yi-epsilon-esti-inball-Bxj}, which, together with \eqref{eq:yi-epsilon-comparison}, Fatou's lemma, the limit in \eqref{eq:yi-epsilon-limit}, and the definition of $y^\epsilon$,  will be employed to prove the desired estimates. 
Indeed for \eqref{eq:exponential-esti-Dirac-many}, by definition of $y^\epsilon$, there holds
\begin{align*}
	\int_{\Omega} \exp[\frac{(4 \pi - \alpha)}{\bm\omega_{\max}}{y}^\epsilon(x)]  dx & = \int_{\Omega} \prod_{i=1}^{k} \exp[\frac{(4 \pi - \alpha)\omega_i}{\bm\omega_{\max}} {y}_i^\epsilon(x)]  dx \\
	& \leq \int_{B_R} \prod_{i=1}^{k} \exp[\frac{(4 \pi - \alpha)\omega_i}{\bm\omega_{\max}} \tilde {y}_i^\epsilon(x)]  dx,
\end{align*}
where the last inequality follows from \eqref{eq:yi-epsilon-comparison} and the fact that $\Omega \subset B_R$.
Setting $B_0 := \cup_{i=1}^k B_{\R^2}(x_i, r_0)$, we can rewrite the right-hand side of the previous estimate as 
\begin{multline*}
	\int_{B_R\backslash B_0} \prod_{i=1}^{k} \exp[\frac{(4 \pi - \alpha)\omega_i}{\bm\omega_{\max}} \tilde {y}_i^\epsilon(x)]  dx \\
	+ \sum_{l=1}^{k} \int_{B_{\R^2}(x_l, r_0)} \prod_{i=1}^{k} \exp[\frac{(4 \pi - \alpha)\omega_i}{\bm\omega_{\max}} \tilde {y}_i^\epsilon(x)]  dx =: A_1^\epsilon + A_2^\epsilon.
\end{multline*}
For $A_1^\epsilon$, we deduce from \eqref{eq:yi-epsilon-esti-pointwise-outside} that
\begin{align*}
	A_1^\epsilon & \leq \int_{B_R\backslash B_0} dx \prod_{i=1}^{k}   ( \frac{2R}{r_0 -\epsilon} )^{(2 - \frac{\alpha}{2 \pi}) \frac{\omega_i}{\bm\omega_{\max}} } 
	%\\
	%&= ( \frac{2R}{r_0 -\epsilon} )^{(2 - \frac{\alpha}{2 \pi}) \sum_{i=1}^{k} \frac{\omega_i}{\bm\omega_{\max}} }  \int_{B_R\backslash B_0} dx \\
	= \pi  ( \frac{2R}{r_0 -\epsilon} )^{(2 - \frac{\alpha}{2 \pi}) \sum_{i=1}^{k} \frac{\omega_i}{\bm\omega_{\max}} } [R^2 -k r_0^2].
\end{align*}
For $A_2^\epsilon$, in view of \eqref{eq:yi-epsilon-esti-pointwise-outside}  and \eqref{eq:yi-epsilon-esti-inball}, there holds
\begin{align*}
	A_2^\epsilon & = \sum_{l=1}^{k} \int_{B_{\R^2}(x_l, r_0)} \exp[\frac{(4 \pi - \alpha)\omega_l}{\bm\omega_{\max}} \tilde {y}_l^\epsilon(x)]  \times   \prod_{i \neq l}  \exp[\frac{(4 \pi - \alpha)\omega_i}{\bm\omega_{\max}} \tilde {y}_i^\epsilon(x)]  dx \\
	& \leq \sum_{l=1}^{k}  ( \frac{2R}{r_0 -\epsilon} )^{(2 - \frac{\alpha}{2 \pi}) \sum_{i \neq l} \frac{\omega_i}{\bm\omega_{\max}} } \int_{B_{\R^2}(x_l, r_0)} \exp[\frac{(4 \pi - \alpha)\omega_l}{\bm\omega_{\max}} \tilde {y}_l^\epsilon(x)]dx \\
	& \leq \sum_{l=1}^{k}  ( \frac{2R}{r_0 -\epsilon} )^{(2 - \frac{\alpha}{2 \pi}) \sum_{i \neq l} \frac{\omega_i}{\bm\omega_{\max}} } \times \frac{2 \pi (\epsilon + r_0)^2}{2 - {(2 - \frac{\alpha}{2 \pi}) \frac{\omega_l}{\bm\omega_{\max}} }} \times ( \frac{2R}{\epsilon + r_0} )^{(2 - \frac{\alpha}{2 \pi}) \frac{\omega_l}{\bm\omega_{\max}} } \\
	& = 2 \pi (r_0 + \epsilon)^2 ( \frac{2R}{r_0 -\epsilon} )^{(2 - \frac{\alpha}{2 \pi})  \sum_{i =1}^k \frac{\omega_i}{\bm\omega_{\max}} } \\
	%\MoveEqLeft[-5] 
	& \qquad \times   \sum_{l=1}^k (\frac{r_0 - \epsilon}{r_0 + \epsilon})^{(2 - \frac{\alpha}{2 \pi}) \frac{\omega_l}{\bm\omega_{\max}} } [{2 - {(2 - \frac{\alpha}{2 \pi}) \frac{\omega_l}{\bm\omega_{\max}} }}]^{-1}.
\end{align*}  
Consequently, one has
\begin{multline*}
	\int_{\Omega} \exp[\frac{(4 \pi - \alpha)}{\bm\omega_{\max}}{y}^\epsilon(x)]  dx \leq
	\int_{B_R} \exp[\frac{(4 \pi - \alpha)}{\bm\omega_{\max}}{y}^\epsilon(x)]  dx \\
	\begin{aligned}
		& \leq \pi  ( \frac{2R}{r_0 -\epsilon} )^{(2 - \frac{\alpha}{2 \pi}) \sum_{i=1}^{k} \frac{\omega_i}{\bm\omega_{\max}} } [R^2 - kr_0^2] 
		+  2 \pi (r_0 + \epsilon)^2 ( \frac{2R}{r_0 -\epsilon} )^{(2 - \frac{\alpha}{2 \pi})  \sum_{i =1}^k \frac{\omega_i}{\bm\omega_{\max}} } \\
		%\MoveEqLeft[-5] 
		&\qquad \times  \sum_{l=1}^k (\frac{r_0 - \epsilon}{r_0 + \epsilon})^{(2 - \frac{\alpha}{2 \pi}) \frac{\omega_l}{\bm\omega_{\max}} } [{2 - {(2 - \frac{\alpha}{2 \pi}) \frac{\omega_l}{\bm\omega_{\max}} }}]^{-1}.
	\end{aligned}
\end{multline*}
This, along with \eqref{eq:yi-epsilon-limit}, the nonnegativity of $y^\epsilon$, and Fatou's lemma, yields \eqref{eq:exponential-esti-Dirac-many}.

It remains to prove \eqref{eq:exponential-esti-Dirac-many-Bxj}. For that purpose, by definition of $y^\epsilon$, there holds
\begin{multline*}
	\int_{B_{\R^2}(x_j, r_0)} \exp[\frac{(4 \pi - \alpha)}{\omega_{j}}{y}^\epsilon(x)]  dx = \int_{B_{\R^2}(x_j, r_0)} \prod_{i=1}^{k} \exp[\frac{(4 \pi - \alpha)\omega_i}{\omega_j} {y}_i^\epsilon(x)]  dx \\
	\begin{aligned}
		& \leq  \int_{B_{\R^2}(x_j, r_0)} \exp[\frac{(4 \pi - \alpha)\omega_j}{\omega_{j}} \tilde {y}_j^\epsilon(x)]  \times   \prod_{i \neq j}  \exp[\frac{(4 \pi - \alpha)\omega_i}{\omega_{j}} \tilde {y}_i^\epsilon(x)]  dx \\
		& \leq  ( \frac{2R}{r_0 -\epsilon} )^{(2 - \frac{\alpha}{2 \pi}) \sum_{i \neq j} \frac{\omega_i}{\omega_{j}} } \times \frac{2 \pi (\epsilon + r_0)^2}{\frac{\alpha}{2 \pi} } \times ( \frac{2R}{\epsilon + r_0} )^{(2 - \frac{\alpha}{2 \pi})  },
	\end{aligned}
\end{multline*}
in view of \eqref{eq:yi-epsilon-limit}, \eqref{eq:yi-epsilon-esti-pointwise-outside-Bxj}  and \eqref{eq:yi-epsilon-esti-inball-Bxj}.  From this, along with \eqref{eq:yi-epsilon-limit}, the nonnegativity of $y^\epsilon$, and Fatou's lemma, we have \eqref{eq:exponential-esti-Dirac-many-Bxj}.
\qed

We now state a direct consequence of \Cref{prop:Poisson-exponential-esti-many-Dirac}, which will contribute significantly to demonstrating the differentiability of the control-to-state operator investigated in \Cref{sec:state-equation}.
\begin{corollary}
	\label{cor:omega-max-4pi}
	Let $g \in L^p(\Omega)$ be arbitrary, but fixed with $p>1$. 	Assume that $\bm\omega \in \R^k$ satisfying  $\bm\omega_{\max} < 4 \pi$. 
	Let $y$ be the unique weak solution to 
	\begin{equation}
		\label{eq:Poisson-many-Dirac-Lp}
		-\Delta y  = g + \sum_{i=1}^k \omega_i \delta_{x_i} \, \text{in } \Omega, \quad
		y  = 0 \, \text{on } \Gamma.
%		
%		\left\{
%		\begin{aligned}
%			-\Delta y & = g + \sum_{i=1}^k \omega_i \delta_{x_i} && \text{in } \Omega, \\
%			y &= 0 && \text{on } \Gamma.
%		\end{aligned}
%		\right.
	\end{equation}
	Then, there exists a constant $C>0$ independent of $\bm{\omega}$ satisfying following statements:
	\begin{enumerate}[label=(\alph*)]
		\item \label{item:omega-max-zero}
		If $\bm\omega_{\max} \leq 0$, then there holds
		\begin{equation}
			\label{eq:exponential-state-belong-Linfty}
			\norm{e^y}_{L^\infty(\Omega)} \leq e^{C\norm{g}_{L^p(\Omega)}}.
		\end{equation}
		\item If $\bm\omega_{\max} \in (0,4\pi)$, then there holds
		\begin{equation}
			\label{eq:exponential-state-belong-L1plus-Poisson-norm}
			\norm{e^y}_{L^{1+\tau}(\Omega)} \leq \pi^{\frac{1}{1+\tau}} e^{C\norm{g}_{L^p(\Omega)}}(\frac{2R}{r_0})^{2k}[R^2 +  \frac{k(1+\tau)\bm{\omega}_{\max} r_0^2}{4\pi -(1+\tau)\bm{\omega}_{\max} }]^{\frac{1}{1+\tau}}
		\end{equation}
		for all $\tau \in [0, \frac{4\pi}{\bm\omega_{\max}} -1 )$, where $R$ and $r_0$ are the constants defined in \Cref{prop:Poisson-exponential-esti-many-Dirac}.
	\end{enumerate}
\end{corollary}
{\it Proof.}  
	Let $y_0$ and $y_1$ be the unique solutions to 
	\begin{align*}
		& 	\left\{
		\begin{aligned}
			-\Delta y_0 & = g  && \text{in } \Omega, \\
			y_0 &= 0 && \text{on } \Gamma
		\end{aligned}
		\right. 
		\quad
		\text{and} \quad
		\left\{
		\begin{aligned}
			-\Delta y_1 & =  \sum_{i=1}^k \omega_i \delta_{x_i} && \text{in } \Omega, \\
			y_1 &= 0 && \text{on } \Gamma,
		\end{aligned}
		\right.
	\end{align*}	
	respectively. 	
	By the linearity, one has $y = y_0 + y_1$. Since $g \in L^p(\Omega)$ with $p> 1 = \frac{N}{2}$, there hold $y_0 \in C(\overline{\Omega})$ and
	\begin{equation}
		\label{eq:y0-Linfty}
		\norm{y_0}_{L^\infty(\Omega)} \leq C \norm{g}_{L^p(\Omega)}
	\end{equation}
	for some constant $C> 0$; see, e.g. \cite[Thm.~4.7]{Troltzsch2010}. 

	We proceed to estimate $e^{y_1}$, and consequently  $e^y$, by considering the following two cases.

	\noindent$\bullet${\it  Case 1:  $\bm\omega_{\max} \leq 0$.}
	In this situation, we deduce from the maximum principle that $y_1 \leq 0$ a.e. in $\Omega$. This, together with \eqref{eq:y0-Linfty}, thus yields \eqref{eq:exponential-state-belong-Linfty}.

	\noindent$\bullet${\it Case 2: $\bm\omega_{\max} \in (0, 4\pi)$.}
	We also derive from the maximum principle that 
	$y_1(x) \leq z(x)$ for a.e. $x \in \Omega$, 
	where $z$ uniquely solves
	\[
		-\Delta z  =  \sum_{i=1}^k \omega_i^+ \delta_{x_i} \quad \text{in } \Omega, \quad 
	z = 0 \quad \text{on } \Gamma.
	\]
	Let $ \tau \in \left[0, \frac{4\pi}{\bm\omega_{\max}} -1 \right)$ be arbitrary.
	Setting now $\alpha : = 4 \pi - (1+ \tau) \bm \omega_{\max} \in (0, 4\pi)$ gives $\frac{(4\pi - \alpha)}{\bm \omega_{\max}}  = 1 + \tau$.
	Applying now \Cref{prop:Poisson-exponential-esti-many-Dirac} to $z$  yields 
	\begin{multline*}
		\int_\Omega e^{(1+\tau)|z|}dx  \\		
		\begin{aligned}
			& \leq  \pi (\frac{2R}{r_0})^{(2 -\frac{\alpha}{2 \pi})\sum_{i=1}^k \frac{\omega_i^{+}}{\bm\omega_{\max} }}[ (R^2 - kr_0^2) + 2r_0^2 \sum_{i=1}^{k} ( 2 - (2 - \frac{\alpha}{2 \pi}) \frac{\omega_i^{+}}{\bm\omega_{\max} } )^{-1} ] \\
			&  \leq \pi (\frac{2R}{r_0})^{k(2 -\frac{\alpha}{2 \pi})}[R^2 + kr_0^2(\frac{4\pi}{\alpha} - 1)] \\
			& = \pi (\frac{2R}{r_0})^{\frac{k(1+\tau)\bm{\omega}_{\max}}{2\pi}}[R^2 +  \frac{k(1+\tau)\bm{\omega}_{\max} r_0^2}{4\pi -(1+\tau)\bm{\omega}_{\max} }],
		\end{aligned}	
	\end{multline*}
	where we have just used the estimates $\omega_i^+ \leq \bm\omega_{\max}  < 4\pi$, $1 \leq i \leq k$, and $( 2 - (2 - \frac{\alpha}{2 \pi}) \frac{\omega_i^{+}}{\bm\omega_{\max} } )^{-1}\leq \frac{2\pi}{\alpha}$
	to obtain the last inequality, and employed the equation $\alpha  = 4 \pi - (1+ \tau) \bm \omega_{\max}$ to derive the last identity.
	Consequently, there holds
	\begin{equation*}
		\int_\Omega e^{(1+\tau)y_1}dx   \leq \pi (\frac{2R}{r_0})^{2k(1+\tau)}[R^2 +  \frac{k(1+\tau)\bm{\omega}_{\max} r_0^2}{4\pi -(1+\tau)\bm{\omega}_{\max} }].
	\end{equation*}
	Combining this with \eqref{eq:y0-Linfty} yields \eqref{eq:exponential-state-belong-L1plus-Poisson-norm}.
\qed

\medskip

%In order to investigate the differentiability of the control-to-state mapping, we need show the existence and uniqueness of solutions to the linearization of the state equation.
%The following result 
%for a linear elliptic equation with measure data and a zeroth-order coefficient belonging to $L^r(\Omega)$ for some $r>1$,
%can be derived from the same argument exploited in \cite[Thm.~4]{AlibertRaymond1997}, where the authors considered the Robin boundary value  linear elliptic partial differential problems with measure data.
%%However, for the convenience, we provide a detailed proof here using Leray--Schauder's principle \cite[Thm.~6.A]{Zeidler_vol1}.
%Its proof is thus skipped.

To investigate the differentiability  of the control-to-state mapping on a certain open set in $\R^k$, we now present the existence and uniqueness of solutions to the linearized state equation. The following result, concerning a linear elliptic equation with measure data and a zeroth-order coefficient in 
$L^r(\Omega)$ for some $r>1$, can be derived using the same argument as in~\cite[Thm.~4]{AlibertRaymond1997}, where the authors treated Robin boundary value problems for linear elliptic PDEs with measure data. 
%For this reason, the proof is omitted.
\begin{lemma}
	\label{lem:linearization-state-equation}
	Assume that $y \in L^1(\Omega)$ satisfying $e^y \in L^r(\Omega)$ with $r >1$. Then, for any $\bm{\eta} \in \R^k$ and $h \in L^1(\Omega)$, there exists a unique $z \in W^{1,q}_0(\Omega)$ for all $1 \leq q < 2$ such that
	\begin{equation}
		\label{eq:linearization-state-equation}
		-\Delta z + e^y z  = h + \sum_{i=1}^k \eta_i \delta_{x_i} \, \text{in } \Omega \quad
		z = 0 \, \text{on } \Gamma.
%		
%		\left\{
%		\begin{aligned}
%			-\Delta z + e^y z & = h + \sum_{i=1}^k \eta_i \delta_{x_i} && \text{in } \Omega \\
%			z &= 0 && \text{on } \Gamma.
%		\end{aligned}
%		\right.
	\end{equation}
	Moreover, there exists a constant $C>0$ independent of $y$, $h$, and $\bm{\eta}$ satisfying
	\begin{equation}
		\label{eq:linearization-state-equation-a-priori-esti}
		\norm{z}_{W^{1,q}_0(\Omega)} \leq C [\norm{h}_{L^1(\Omega)} + \sum_{i=1}^k |\eta_i| ].
	\end{equation}
\end{lemma}
{\it Proof.}  
	Let  $r'$ be the conjugate number of $r$, i.e., $r' : = \frac{r}{r-1}$ and let $q \in [1, 2)$ be close to $2$, such that $r' < \frac{2q}{2-q}$.	
	The choice of $q$ ensures that $W^{1,q}_0(\Omega)$ is compactly embedded in $L^{r'}(\Omega)$, due to the Rellich--Kondrachov theorem; see, e.g. \cite[Thm.~6.2]{Adams1975}.
	We consider the mapping 
	\[
	\begin{aligned}
		T:  L^{r'}(\Omega) & \to W^{1,q}_0(\Omega) \hookrightarrow L^{r'}(\Omega) \\
		v & \mapsto T(v) := z_v,
	\end{aligned}
	\]
	where $z_v$ is the unique solution to
	\begin{equation}
		\label{eq:-linearization-state-equation-Tv}
		-\Delta z_v + e^y v  = h +  \sum_{i=1}^k \eta_i \delta_{x_i} \, \text{in } \Omega \quad
		z_v = 0 \, \text{on } \Gamma.
%		
%		\left\{
%		\begin{aligned}
%			-\Delta z_v + e^y v & = h +  \sum_{i=1}^k \eta_i \delta_{x_i} && \text{in } \Omega \\
%			z_v &= 0 && \text{on } \Gamma.
%		\end{aligned}
%		\right.
	\end{equation}
	We first show that $T$ is continuous as a mapping from $L^{r'}(\Omega)$ to itself. Indeed, for any $v_1, v_2 \in L^{r'}(\Omega)$, by setting $z_i := T(v_i)$ and then subtracting the equations for $z_1$ and $z_2$, there holds
	\[
	-\Delta (z_1 - z_2) = -e^{y}(v_1 - v_2) \quad \text{in } \Omega, \quad z_1 - z_2 =0 \quad \text{on } \Gamma.
	\]
	Obviously, the right-hand side of this equation belongs to $L^1(\Omega)$. Using the a priori estimate for Poisson's equation with $L^1$-data yields
	\[
	\norm{z_1 -z_2}_{W^{1,q}_0(\Omega)} \leq C\norm{e^y(v_1 - v_2)}_{L^1(\Omega)} \leq C \norm{e^y}_{L^r(\Omega)}\norm{v_1-v_2}_{L^{r'}(\Omega)};
	\]
	see \cite{BrezisStrauss1973}, where we have exploited the H\"{o}lder inequality to derive the last estimate. From this and the continuous imbedding $W^{1,q}_0(\Omega) \hookrightarrow L^{r'}(\Omega)$, the continuity of $T$ then follows. 
	
	Since $W^{1,q}_0$ is compactly embedded in $L^{r'}(\Omega)$, $T: L^{r'}(\Omega) \to L^{r'}(\Omega)$ is compact. 
	
	We now show that the following set
	\[
	K := \{ z \in L^{r'}(\Omega) \, | \, \exists t \in (0,1): z = T(tz) \}
	\]
	is bounded. To this end, for any $z \in K$, there holds 
	\begin{equation}
		\label{eq:K-bounded-pdes}
		-\Delta z + te^y z  = h+ \sum_{i=1}^k \eta_i \delta_{x_i} \,\text{in } \Omega \quad
		z = 0 \, \text{on } \Gamma
%		
%		\left\{
%		\begin{aligned}
%			-\Delta z + te^y z & = h+ \sum_{i=1}^k \eta_i \delta_{x_i} && \text{in } \Omega \\
%			z &= 0 && \text{on } \Gamma
%		\end{aligned}
%		\right.
	\end{equation}
	for some $t \in [0,1]$. By \cite[Prop.~4.B.3]{BresisMarcusPonce2007} for $f := te^{y}z$, one has
	\[
	\norm{te^y z }_{L^1(\Omega)} = \int_\Omega te^y z \sign(z) dx \leq \norm{h}_{L^1(\Omega)} +  \sum_{i=1}^k |\eta_i|.
	\]
	The a priori estimate of solutions to \eqref{eq:K-bounded-pdes}; see, e.g. \cite[Thm.~4.B.1]{BresisMarcusPonce2007}, gives
	\[
	\norm{z}_{W^{1,q}_0(\Omega)} \leq C [\norm{te^y z}_{L^1(\Omega)} + \norm{h}_{L^1(\Omega)} + \sum_{i=1}^k |\eta_i| ] \leq C [\norm{h}_{L^1(\Omega)} +\sum_{i=1}^k |\eta_i|],
	\]
	which, along with $W^{1,q}_0(\Omega) \hookrightarrow L^{r'}(\Omega)$, thus shows the boundedness of $K$.
	
	\medskip
	
	We now apply the Leray--Schauder principle to mapping $T$ and deduce the existence of solutions to \eqref{eq:linearization-state-equation}. 
	On the other hand, the uniqueness of solutions to \eqref{eq:linearization-state-equation} follows from \cite[Cor.~ 4.B.1]{BresisMarcusPonce2007}.
	Finally, by \cite[Prop.~4.B.3]{BresisMarcusPonce2007}, there holds
	\begin{align*}
		\norm{e^y z }_{L^1(\Omega)} = \int_{\Omega} e^y z \sign z dx  \leq \norm{h}_{L^1(\Omega)} +  \sum_{i=1}^k |\eta_i|.
	\end{align*}
	From this, we can conclude from \cite[Thm.~4.B.1]{BresisMarcusPonce2007} that
	\begin{align*}
		\norm{z}_{W^{1,q}_0(\Omega)}  & \leq C [\norm{e^y z}_{L^1(\Omega)} + \norm{h}_{L^1(\Omega)} + \sum_{i=1}^k |\eta_i| ] 
		%\\
		%& 
		\leq C [\norm{h}_{L^1(\Omega)} +\sum_{i=1}^k |\eta_i|],
	\end{align*}
	which gives \eqref{eq:linearization-state-equation-a-priori-esti}.
\qed

\medskip 

%%% The continuity of Nemytskii's operator induced by the exponential function
With the aim of studying the Fr\'{e}chet differentiability of the control-to-state operator, considered in \Cref{sec:state-equation}, we need to introduce some open sets,  depending on $\bm{b} = (b_1, b_2, \ldots, b_k) \in \R^k$ with $\bm{b}_{\max} \in (0, 4\pi)$, and Banach spaces as follows: 
%Let us consider the open subsets, depending on $\bm{b} = (b_1, b_2, \ldots, b_k) \in \R^k$ with $\bm{b}_{\max} \in (0, 4\pi)$, respectively,  in $\R^k$, $\mathcal{M}(\Omega)$, and $W^{1,1}_0(\Omega)$ as follows
\begin{align}
	\label{eq:Ob-set}
	\mathcal{O}_{\bm{b}} & :=   \{\bm\omega = (\omega_1, \omega_2, \ldots, \omega_k) \in \R^k \, | \, \omega_i < b_i \, \forall 1 \leq i \leq k  \}, \\
	\mathcal{D}_{\bm{b}} & :=  \{ \mu = \sum_{i =1}^k \omega_i \delta_{x_i} \in \mathcal{M}(\Omega) \, | \, \bm\omega = (\omega_1, \omega_2, \ldots, \omega_k) \in \mathcal{O}_{\bm{b}} \}, \label{eq:Db-set} \\
	V_{\bm{b}} & :=  \{ y \in W^{1,1}_0(\Omega) \,| \, - \Delta y \in \mathcal{D}_{\bm{b}}  + L^{r_{\bm{b}}}(\Omega)  \},  \label{eq:Vb-set} 
	\intertext{and}
	\mathcal{D} & := \{ \mu = \sum_{i =1}^k \omega_i \delta_{x_i} \in \mathcal{M}(\Omega) \,| \, \bm\omega = (\omega_1, \omega_2, \ldots, \omega_k) \in \R^k  \}, \label{eq:D-space} \\
	W_{r_{\bm{b}}} & := \{ y \in W^{1,1}_0(\Omega) \, | \, - \Delta y \in \mathcal{D}  + L^{r_{\bm{b}}}(\Omega)   \} \label{eq:Wrb-space}
\end{align}
with
\begin{equation}
	\label{eq:rb-constant}
	r_{\bm{b}}:= \min \{\frac{4\pi}{{\bm{b_{\max}}}}, p \}  >1.
\end{equation}
Obviously, $\mathcal{D}$ is a closed finite-dimensional linear subspace of $\mathcal{M}(\Omega)$, while $W_{r_{\bm{b}}}$ is a Banach space when endowed with the norm
\[
\norm{y}_{W_{r_{\bm{b}}}} := \norm{y}_{W^{1,1}_0(\Omega)} + \norm{\Delta y}_{\mathcal{D}  + L^{r_{\bm{b}}}(\Omega)}.
\]
Here, for any $y \in W_{r_{\bm{b}}}$, the $(\mathcal{D}  + L^{r_{\bm{b}}}(\Omega))$-norm of $\Delta y$ is given by
\[
\norm{\Delta y}_{\mathcal{D}  + L^{r_{\bm{b}}}(\Omega)} = \norm{\mu}_{\mathcal{M}(\Omega)} + \norm{h}_{L^{r_{\bm{b}}}(\Omega)},
\] 
where, $(\mu, h) \in (\mathcal{D}, L^{r_{\bm{b}}}(\Omega))$ is the unique pair satisfying $-\Delta y = \mu + h$. 
We also observe from \Cref{lem:exp-V-Omega} below that the space $W_{r_{\bm{b}}}$ is continuously embedded in $W^{1,q}_0(\Omega)$ for all $1 \leq q < \frac{N}{N-1} = 2$. 
Furthermore, by \Cref{lem:openness-Db-wWrb},
$\mathcal{D}_{\bm{b}}$ and 
$V_{\bm{b}}$ are, respectively, open subsets of $\mathcal{D}$ and $W_{r_{\bm{b}}}$.

\medskip 
%%% The continuity of The Nemytskii's operator associated with the exponential function defined over $V_{\bm{ b }}$
The continuity of  Nemytskii's operator associated with the exponential function as a mapping defined on $V_{\bm{ b }}$ is established below and plays a crucial role in displaying the Fr\'{e}chet differentiability of the control-to-state mapping.
\begin{proposition}
	\label{prop:Nemytskii-exp-continuity}
	Assume that $\bm{b}_{\max} \in (0, 4\pi)$ and let  $y_* \in V_{\bm{b}}$ be arbitrary.
	Then, there exist positive constants $\rho_* = \rho_*(y_*)$ and $\epsilon_* = \epsilon_*(y_*)$ such that
	\begin{equation}
		\label{eq:exp-y*-Lr}
		e^{y} \in {L^{r_{\bm{ b }} +\epsilon_*  }(\Omega)} \quad \text{for all } y \in W_{r_{\bm{b}}} \quad \text{such that } \norm{y - y_*}_{W_{r_{\bm{b}}}} <  \rho_*.
	\end{equation} 
	Furthermore, there holds
	\begin{equation}
		\label{eq:Nemy-exp-continuity-Vb}
		\norm{e^y - e^{y_*}}_{L^{r_{\bm{ b }} +\epsilon_*  }(\Omega)} \to 0 \quad \text{as} \quad \norm{y - y_*}_{W_{r_{\bm{b}}}} \to 0.
	\end{equation}
	% The Nemytskii's operator associated with the exponential function is continuous as a mapping from $V_{\bm{ b }}$.
\end{proposition}
{\it Proof.}  
	By fixing now $y_* \in V_{\bm{b}}$, there exists a unique pair $(\bm{\omega}^*, g_*) \in \R^k \times L^{r_{\bm{ b }}}(\Omega)$ with $\bm{\omega}^{*} = (\omega^*_1, \omega_2^*, \ldots, \omega_k^*) \in \R^k$ such that
	\begin{align*}
		& \omega_i^* < b_i \quad \text{for all } 1 \leq i \leq k 
		\quad \text{and} \quad -\Delta y_* = g_* + \sum_{i=1}^k \omega_i^* \delta_{x_i} \quad \text{in } \Omega.
%		\intertext{and}
%		& -\Delta y_* = g_* + \sum_{i=1}^k \omega_i^* \delta_{x_i} \quad \text{in } \Omega.
	\end{align*}
	Similarly, for any $y \in W_{r_{\bm{b}}}$, by definition, there exists a unique pair  $(\bm{\omega}, g) \in \R^k \times L^{r_{\bm{ b }}}(\Omega)$ satisfying
	$-\Delta y = g + \sum_{i=1}^k \omega_i \delta_{x_i}$ in $\Omega$.
	
%	\[
%	-\Delta y = g + \sum_{i=1}^k \omega_i \delta_{x_i} \quad \text{in } \Omega.
%	\]
	By definition, there thus holds
	\[
	\norm{y - y_*}_{W_{r_{\bm{b}}}} = \norm{y - y_*}_{W^{1,1}_0(\Omega)} + \norm{g - g_*}_{L^{r_{\bm{ b }}}(\Omega)} + \sum_{i =1}^k | \omega_i - \omega_i^*|.
	\]
	Fix $\delta_* >0$ such that $0 < \bm{\omega}_{\max} +  \delta_*  < \bm{ b }_{\max}$
	and choose  
	\begin{equation*}
		%\label{eq:epsilon-star}
		\epsilon_* := \frac{1}{2} [\frac{4\pi}{\bm{\omega}^*_{\max} +  \delta_* } - \frac{4 \pi}{\bm{ b }_{\max}}] > 0.
	\end{equation*} 
	Easily,  \eqref{eq:rb-constant} implies
	$ 
	2\epsilon_* + r_{\bm{ b }} \leq  \frac{4\pi}{\bm{\omega}^*_{\max} +  \delta_* }.
	$		
	Assume now that $\norm{y - y_*}_{W_{r_{\bm{b}}}} < \frac{\delta_*}{2}$. This gives 
	$
	\sum_{i =1}^k | \omega_i - \omega_i^*| < \frac{\delta_*}{2},
	$
	and consequently, 
	\[
	\omega_i < \frac{\delta_*}{2} + \omega_i^* \leq \frac{\delta_*}{2} + \bm{\omega}^*_{\max}  \quad \text{for all } 1 \leq i \leq k.
	\]
	We therefore have
	\[
	(2\epsilon_* + r_{\bm{ b }})\bm{\omega}_{\max} \leq \frac{4\pi}{\bm{\omega}^*_{\max} +  \delta_* }(\frac{\delta_*}{2} + \bm{\omega}^*_{\max} ) < 4 \pi 
	\]
	and thus
	\[
		4\pi - (2\epsilon_* + r_{\bm{ b }})\bm{\omega}_{\max} \geq 4\pi - \frac{4\pi}{\bm{\omega}^*_{\max} +  \delta_* }(\frac{\delta_*}{2} + \bm{\omega}^*_{\max} ) > 0.
	\]
	Applying now \eqref{eq:exponential-state-belong-L1plus-Poisson-norm} for $\tau := (2\epsilon_* + r_{\bm{ b }} -1)$ and $p := r_{\bm{ b }}$ yields that 
	\begin{equation}
		\label{eq:exp-y-bounded-Lplus}
		\norm{e^y}_{L^{r_{\bm{ b }} + 2\epsilon_*}(\Omega)} \leq C(R,r_0) e^{C\norm{g}_{ L^{r_{\bm{ b }} }(\Omega)}}\leq  C(R,r_0)e^{C(\norm{g_*}_{L^{r_{\bm{ b }}}(\Omega) } + \frac{\delta_*}{2})}
	\end{equation}
	for all $y \in W_{r_{\bm{b}}}$ such that $\norm{y - y_*}_{W_{r_{\bm{b}}}} < \frac{\delta_*}{2}$.
	Putting $\rho_* := \frac{\delta_*}{2}$ and exploiting the continuous embedding $L^{r_{\bm{ b }} + 2\epsilon_*}(\Omega) \hookrightarrow L^{r_{\bm{ b }} + \epsilon_*}(\Omega)$ yield \eqref{eq:exp-y*-Lr}.

	%We have shown that $e^y \in L^{r_{\bm{ b }} + \epsilon_*}(\Omega)$ whenever $\norm{y - y_*}_{W_{r_{\bm{b}}}} < \rho_*$. 
	It remains to show \eqref{eq:Nemy-exp-continuity-Vb}. For this purpose, let $\norm{y - y_*}_{W_{r_{\bm{b}}}} \to 0$. 
	From the definition of the norm in $W_{r_{\bm{b}}}$ (see \eqref{eq:Wrb-space}), there thus holds $y \to y_*$  and thus $e^{y} \to e^{y_*}$ a.e. in $\Omega$. 
	Consequently, $e^{y} \to e^{y_*}$ in measure in $\Omega$.
	From this and \eqref{eq:exp-y-bounded-Lplus}, we obtain \eqref{eq:Nemy-exp-continuity-Vb} from \cite[Chap.~XI, Prop.~3.10]{Visintin1996}.
\qed   

The G\^{a}teaux differentiability at the origin of the Nemytskii operator, associated with the exponential function over $W_{r_{\bm{b}}}$, is partly stated below.
\begin{lemma}
	\label{lem:exp-V-Omega}
	For any $z \in W_{r_{\bm{b}}}$ and any $q \in [1,2)$, $r \geq 1$, there hold
	\begin{align}
		\norm{z}_{W^{1,q}_0(\Omega)}  \leq C_1 \norm{z}_{W_{r_{\bm{b}}}}, \quad  \norm{z}_{L^r(\Omega)} \leq C_2 \norm{z}_{W_{r_{\bm{b}}}} \label{eq:W1q-Wrb-embedding-V-Omega-Lr-esti}
		%\norm{z}_{L^r(\Omega)} & \leq C_2 \norm{z}_{W_{r_{\bm{b}}}} \label{eq:V-Omega-Lr-esti}
	\end{align}
	for some constants $C_1, C_2>0$,
	and
	\begin{equation}
		\label{eq:exp-Nemy-G-diff-at-origin}
		\norm{\frac{e^{tz} -1}{t}  - z}_{L^{r}(\Omega)} \to 0 \quad \text{as} \quad t \to 0^+.
	\end{equation}
\end{lemma}
{\it Proof.}  
	Let $z \in W_{r_{\bm{b}}}$ and $r \geq 1$ be arbitrary, but fixed. By definition, there exists a pair $(\bm{\omega}, g) \in \R^k \times L^{r_{\bm{ b }}}(\Omega)$ satisfying
	\[
	-\Delta z = g + \sum_{i=1}^k \omega_i \delta_{x_i} \quad \text{in } \Omega.
	\]
	From \cite[Thm.~4.B.1]{BresisMarcusPonce2007}, one has $z \in W^{1,q}_0(\Omega)$ for all $1 \leq q < 2$
	and there holds
	\[
	\norm{z}_{W^{1,q}_0(\Omega)} \leq C \norm{g + \sum_{i=1}^k \omega_i \delta_{x_i}}_{\mathcal{M}(\Omega)} \leq C ( \norm{g}_{L^1(\Omega)} + \sum_{i =1}^k |\omega_i| ),
	\]
	which, as well as the continuous embedding $L^{r_{\bm{b}}}(\Omega) \hookrightarrow L^1(\Omega)$, gives 
	\[
	\norm{z}_{W^{1,q}_0(\Omega)} \leq C ( \norm{g}_{L^{r_{\bm{b}}}(\Omega)} + \sum_{i =1}^k |\omega_i| ) \leq C \norm{z}_{W_{r_{\bm{b}}}}.
	\]
	Here, we have just used the definition of norm on $W_{r_{\bm{b}}}$ in order to get the last estimate.	
	We therefore obtain the first estimate in \eqref{eq:W1q-Wrb-embedding-V-Omega-Lr-esti}.
	
	Moreover, the first estimate in \eqref{eq:W1q-Wrb-embedding-V-Omega-Lr-esti}, in combination with the continuous embedding $W^{1,q}_0(\Omega) \hookrightarrow L^{r}(\Omega)$ for  $q$ sufficiently close to $2$, then yields the second estimate in \eqref{eq:W1q-Wrb-embedding-V-Omega-Lr-esti}.
	
	It remains to prove \eqref{eq:exp-Nemy-G-diff-at-origin}. To this end, let $t_0 >0$ be fixed such that $t_0 r \bm{\omega}_{\max} < 4 \pi$. Since $t_0z$ satisfies
	\[
	-\Delta (t_0z) = t_0g + \sum_{i=1}^k t_0\omega_i \delta_{x_i} \quad \text{in } \Omega, \quad t_0z = 0 \quad \text{on } \Gamma,
	\]
	there holds 
	\begin{equation}
		\label{eq:exp-z-t0r}
		e^{t_0z} \in L^r(\Omega),
	\end{equation}
	in light of \eqref{eq:exponential-state-belong-Linfty} and \eqref{eq:exponential-state-belong-L1plus-Poisson-norm}. 
	We now define the function depending on $t \in (0, t_0]$,
	\[
	\gamma_t(x) :=  \frac{e^{tz(x)}-1}{t} - z(x)    \quad \text{for a.e. } x \in \Omega.
	\]
	Obviously, one has
	\begin{equation}
		\label{eq:gamma-t-func-limit}
		\gamma_t(x) \to 0 \quad \text{as} \quad t \to 0^+ \quad \text{for a.e. } x \in \Omega.
	\end{equation}
	Since $\psi(t) := \frac{e^{at}-1}{t}$, $a \in \R$, is nondecreasing on $(0, t_0]$, it follows that 
	\[
	0 \leq \gamma_t(x) \leq \frac{e^{t_0z(x)}-1}{t_0} - z(x) \quad \text{for all } t \in (0, t_0] \quad \text{and for a.e. } x \in \Omega.
	\]
	From this, \eqref{eq:gamma-t-func-limit}, and \eqref{eq:exp-z-t0r}, we then deduce from the Lebesgue dominated convergence theorem that $\gamma_t \to 0$ in $L^r(\Omega)$ as $t \to 0^+$. 
	The limit \eqref{eq:exp-Nemy-G-diff-at-origin}  follows.
\qed   

\medskip 
From now on, let us define the mapping 
\begin{equation}
	\label{eq:F0-operator}
	\begin{aligned}
		F_0: V_{\bm{ b }} & \to  \mathcal{D} + L^{r_{\bm{b}}}(\Omega)\\
		y	& \mapsto -\Delta y + e^y - 1,
	\end{aligned}
\end{equation}
where $V_{\bm{ b }}$, $\mathcal{D}$, and $r_{\bm{ b }}$ are determined in \eqref{eq:Vb-set}, \eqref{eq:D-space}, and \eqref{eq:rb-constant}, respectively.

\begin{lemma}
	\label{lem:F0-differentiability}
	Assume that $\bm{b}_{\max} \in (0, 4\pi)$. Then $F_0$ is of class $C^1$ on $V_{\bm{ b }}$. Moreover, its derivative is  given by
	\begin{equation}
		\frac{\partial F_0}{\partial y}(y)z = - \Delta z + e^yz \label{eq:F0-1st-der}
	\end{equation}
	for all $y \in V_{\bm{ b }}$ and for all $z \in W_{r_{\bm{b}}}$.
	%	\begin{align}
		%		& \frac{\partial F_0}{\partial y}(y)z = - \Delta z + e^yz \label{eq:F0-1st-der}
		%		\intertext{and}
		%		&  \frac{\partial^2 F_0}{\partial y^2}(y)[z_1,z_2] =  e^yz_1z_2 \label{eq:F0-2nd-der}
		%	\end{align}
	%	for all $z, z_1, z_2 \in W_{r_{\bm{b}}}$. 
\end{lemma}
{\it Proof.}  
	To prove $F_0$ belongs to the class $C^1$  on  $V_{\bm{ b }}$, we first indicate that $F_0$ is G\^{a}teaux differentiable and then show that its G\^{a}teaux derivative is continuous; see, e.g., \cite[Prop.~2.51]{Penot2013}. We thus split the proof into two parts.
	
	%	%\begin{enumerate}[label=\alph*]
	%		%\item \label{item:1st-F0} 
	%	\noindent$\star$ {\it Step 1:}	\emph{Showing that $F_0$ belongs to $C^1$ class on  $V_{\bm{ b }}$.} This part will be done by indicating that $F_0$ is G\^{a}teaux differentiable and its G\^{a}teaux derivative is continuous and then by applying, e.g., \cite[Prop.~2.51]{Penot2013}. 
	
	%\begin{enumerate}[label= \ref{item:1st-F0}\arabic*] 
	%\item \label{item:F0-G-diff} 
	
	%\noindent$\bullet$ {\it Step 1.1:}
	\noindent$\bullet$ {\it Step 1:}
	\emph{G\^{a}teaux differentiability of $F_0$}. Indeed, for any $y_* \in V_{\bm{ b }}$ and $z \in W_{r_{\bm{b}}}$, there holds for all $t>0$ that
	\begin{multline*}
		\norm{\frac{F_0(y_* + tz) - F_0(y_*)}{t} -  \frac{\partial F_0}{\partial y}(y_*)z}_{ \mathcal{D} + L^{r_{\bm{b}}}(\Omega)}  = \norm{ e^{y_*}( \frac{e^{tz} -1}{t}  - z)}_{ \mathcal{D} + L^{r_{\bm{b}}}(\Omega)} \\
		\begin{aligned}
			& = \norm{ e^{y_*}( \frac{e^{tz} -1}{t}  - z)}_{L^{r_{\bm{b}}}(\Omega)} 
			\leq \norm{ e^{y_*}}_{L^{s}(\Omega)}\norm{\frac{e^{tz} -1}{t}  - z}_{L^{s_1}(\Omega)},
		\end{aligned}
	\end{multline*}
	where 
	\begin{equation}
		\label{eq:s-s1-constant}
		s := r_{\bm{b}} + \epsilon_* \quad \text{and} \quad \frac{1}{s_1} := \frac{1}{r_{\bm{b}}} - \frac{1}{s}
	\end{equation}
	with $\epsilon_*$ being a constant determined in \Cref{prop:Nemytskii-exp-continuity}.
	Letting now  $t \to 0^+$ and exploiting the limit \eqref{eq:exp-Nemy-G-diff-at-origin} in \Cref{lem:exp-V-Omega}, we arrive at 
	\begin{equation}
		\label{eq:F0-G-der-limit}
		\norm{\frac{F_0(y_* + tz) - F_0(y_*)}{t} -  \frac{\partial F_0}{\partial y}(y_*)z}_{ \mathcal{D} + L^{r_{\bm{b}}}(\Omega)}  \to 0.
	\end{equation}
	On the other hand, for any $z \in W_{r_{\bm{b}}}$, we deduce from H\"{o}lder's inequality and the definition of norm in $W_{r_{\bm{b}}}$ that
	\begin{align*}
		\norm{\frac{\partial F_0}{\partial y}(y_*)z}_{\mathcal{D}+ L^{r_{\bm{b}}}(\Omega)} & = \norm{ -\Delta z + e^{y_*}z}_{\mathcal{D}+ L^{r_{\bm{b}}}(\Omega)} 
		\leq  \norm{ -\Delta z}_{\mathcal{D}+ L^{r_{\bm{b}}}(\Omega)}  +  \norm{  e^{y_*}z}_{ L^{r_{\bm{b}}}(\Omega)} \\
		& \leq  \norm{ -\Delta z}_{\mathcal{D}+ L^{r_{\bm{b}}}(\Omega)}  + \norm{e^{y_*}}_{L^s(\Omega)}  \norm{  z}_{ L^{s_1}(\Omega)} \\
		& \leq \norm{  z}_{W_{r_{\bm{b}}}} + \norm{e^{y_*}}_{L^s(\Omega)} \norm{  z}_{W_{r_{\bm{b}}}},
	\end{align*}
	where we have just exploited the second estimate in \eqref{eq:W1q-Wrb-embedding-V-Omega-Lr-esti} to get the last inequality. 
	This confirms that $\frac{\partial F_0}{\partial y}(y_*): W_{r_{\bm{b}}} \to \mathcal{D}+ L^{r_{\bm{b}}}(\Omega)$ is continuous. From this and \eqref{eq:F0-G-der-limit}, the G\^{a}teaux differentiability of $F_0$ at $y_*$ follows.
	Since $y_*$ is arbitrary in $V_{\bm{b}}$, $F_0$ is G\^{a}teaux differentiable over $V_{\bm{b}}$.

	%\noindent$\bullet$ {\it Step 1.2:} 
	\noindent$\bullet$ {\it Step 2:} 
	\emph{Continuity of $\frac{\partial F_0}{\partial y}$}. 
	In fact, we have for any $y, y_* \in V_{\bm{ b }}$ and for constants $s$, $s_1$ defined as in \eqref{eq:s-s1-constant} that
	\begin{multline*}
		\norm{  \frac{\partial F_0}{\partial y}(y) -  \frac{\partial F_0}{\partial y}(y_*)}_{ \Linop( W_{r_{\bm{b}}},\mathcal{D} + L^{r_{\bm{b}}}(\Omega) ) } \\
		\begin{aligned}
			& = \sup\left \{\norm{  \frac{\partial F_0}{\partial y}(y)z -  \frac{\partial F_0}{\partial y}(y_*)z }_{\mathcal{D} + L^{r_{\bm{b}}}(\Omega) } | \norm{z}_{W_{r_{\bm{b}}}} \leq 1 \right \} \\
			&  = \sup\{\norm{ (e^{y}- e^{y_*} )z }_{\mathcal{D} + L^{r_{\bm{b}}}(\Omega) } | \norm{z}_{W_{r_{\bm{b}}}} \leq 1 \} \\
			&  = \sup\{\norm{ (e^{y}- e^{y_*} )z }_{ L^{r_{\bm{b}}}(\Omega) } | \norm{z}_{W_{r_{\bm{b}}}} \leq 1 \} \\
			& \leq \sup\{\norm{ e^{y}- e^{y_*}  }_{ L^{s}(\Omega) } \norm{z}_{L^{s_1}(\Omega)} | \norm{z}_{W_{r_{\bm{b}}}} \leq 1 \} \\
			& \leq C  \norm{ e^{y}- e^{y_*}  }_{ L^{s}(\Omega) },
		\end{aligned}
	\end{multline*}
	where we have just employed the second estimate in \eqref{eq:W1q-Wrb-embedding-V-Omega-Lr-esti}  to get the last estimate. 
	This, along with \Cref{prop:Nemytskii-exp-continuity}, yields the continuity of $\frac{\partial F_0}{\partial y}$ in $y_*$. Since $y_*$ is arbitrary, $\frac{\partial F_0}{\partial y}$ is continuous on $V_{\bm{ b }}$. 
\qed   

\begin{lemma}
	\label{lem:F0-der-isomorphism}
	Under the hypothesis of \Cref{lem:F0-differentiability}, for any $y \in V_{\bm{b}}$, the mapping $\frac{\partial F_0}{\partial y}(y)$ is an isomorphism as an operator from $W_{r_{\bm{b}}}$ to $\mathcal{D} + L^{r_{\bm{ b }}}(\Omega)$. 
\end{lemma}
{\it Proof.}  
	Let $y \in V_{\bm{b}}$ be arbitrary, but fixed.
	By  \eqref{eq:exp-y*-Lr} in  \Cref{prop:Nemytskii-exp-continuity}, there holds
	\[
	e^{y} \in L^r(\Omega)
	\]
	for some $r > r_{\bm{b}}$. 
	Let us define the operator $L: \mathcal{D} + L^{r_{\bm{ b }}}(\Omega) \to W_{r_{\bm{b}}}$ as follows: For any $\mu := \sum_{i=1}^k \omega_i \delta_{x_i} + h$ with $\bm{\omega} = (\omega_1, \ldots, \omega_k) \in \R^k$ and $h \in L^{r_{\bm{b}}}(\Omega)$, $L\mu := z$ is the unique solution in $W^{1,q}_0(\Omega)$, $1 \leq q <2$, to
	\[
	-\Delta z + e^{y}z = \mu \quad \text{in } \Omega;
	\]
	see, \Cref{lem:linearization-state-equation}. 
	The well-posedness of $L$ is as follows.
	%By \Cref{lem:linearization-state-equation}, $z$ is unique and belongs to $W^{1,q}_0(\Omega)$ for all $1 \leq q <2$. 
	Since $W^{1,q}_0(\Omega) \hookrightarrow L^{r_1}(\Omega)$ with $r_1 := (r_{\bm{b}}^{-1} - r^{-1})^{-1}$ and $q$ sufficiently close to $2$, $z \in L^{r_1}(\Omega)$. Then, one has $e^{y}z \in L^{r_{\bm{b}}}(\Omega)$ and thus $z \in W_{r_{\bm{b}}}$.

	It suffices to show the continuity of $T$. To this end, we deduce from \eqref{eq:linearization-state-equation-a-priori-esti} and the continuous embeddings $W^{1,q}_0(\Omega)\hookrightarrow  W^{1,1}_0(\Omega)$ for $1 \leq q < 2$  and $L^{r_{\bm{b}}}(\Omega) \hookrightarrow L^1(\Omega)$ that
	\begin{multline}
		\label{eq:z-W11-esti}
		\norm{z}_{W^{1,1}_0(\Omega)} \leq C_1 \norm{z}_{W^{1,q}_0(\Omega)} \leq C_2 \norm{\mu}_{\mathcal{M}(\Omega)} = C_2 [\norm{h}_{L^{1}(\Omega)}  + \sum_{i =1}^k |\omega_i| ] \\
		\leq C_3 [\norm{h}_{L^{r_{\bm{b}}}(\Omega)}  + \sum_{i =1}^k |\omega_i| ]
	\end{multline}
	for some constants $C_1, C_2, C_3 >0$.
	From this,  the triangle inequality,  H\"{o}lder's inequality, 
	and the embedding $W^{1,q}(\Omega) \hookrightarrow L^{r_1}(\Omega)$ for $q$ sufficiently close to $2$,
	we arrive at
	\begin{align*}
		\norm{\Delta z}_{\mathcal{D} + L^{r_{\bm{ b }}}(\Omega)} & = \norm{h - e^y z + \sum_{i =1}^k \omega_i \delta_{x_i}}_{\mathcal{D} + L^{r_{\bm{ b }}}(\Omega)} %\\
		%& 
		=  \norm{h- e^yz}_{L^{r_{\bm{b}}}(\Omega)}  + \sum_{i =1}^k |\omega_i| \\
		& \leq \norm{h}_{L^{r_{\bm{b}}}(\Omega)} + \norm{e^y}_{L^{r}(\Omega)}\norm{z}_{L^{r_1}(\Omega)}  + \sum_{i =1}^k |\omega_i|  \\
		& \leq \norm{h}_{L^{r_{\bm{b}}}(\Omega)}  + \sum_{i =1}^k |\omega_i| + C \norm{e^y}_{L^r(\Omega)} \norm{z}_{W^{1,q}_0(\Omega)} \\
		& \leq C_4 (1+  \norm{e^y}_{L^r(\Omega)})[\norm{h}_{L^{r_{\bm{b}}}(\Omega)}  + \sum_{i =1}^k |\omega_i| ].
	\end{align*}
	Combining this with \eqref{eq:z-W11-esti} and the definition of the $W_{r_{\bm{b}}}$-norm, we derive the continuity of $L$.
\qed   

\section{State equation}
\label{sec:state-equation}
%	In order to further investigate \eqref{eq:P-original} and \eqref{eq:P-M-prob}, we need to take into account the state equation \eqref{eq:P-state}.
%	\begin{equation}
	%		\tag{\ref{eq:P-state}}
	%		\left\{
	%		\begin{aligned}
		%			-\Delta y + \left( e^y -1  \right)  & = f_0 + \sum_{i=1}^k \eta_i \delta_{x_i} && \text{in } \Omega \\
		%			y &= 0 && \text{on } \Gamma,
		%		\end{aligned}
	%		\right.
	%	\end{equation}
%	with $f_0 \in L^p(\Omega)$, $p >1$.  

Recall from \Cref{prop:state-equation} that  the state equation \eqref{eq:P-state} admits at least a weak solution $y$ if and only if $\bm{\eta}$ satisfies
\eqref{eq:SE-existence-solution-assumption-0}.
%	\begin{equation}
	%		\label{eq:SE-existence-solution-assumption}
	%		\bm\eta_{\max} \leq 4\pi.
	%	\end{equation} 
Moreover, under \eqref{eq:SE-existence-solution-assumption-0}, the solutions to \eqref{eq:P-state} are unique and belong to $W^{1,q}_0(\Omega)$ for any $1 \leq q < \frac{N}{N-1} =2$.
From this, we now define the control-to-state mapping, denoted by $S$, as follows
\begin{equation}
	\label{eq:control-2-state-oper}
	\begin{aligned}
		S: (-\infty, 4 \pi]^k &\to   W^{1,q}_0(\Omega)\\
		\bm{\eta} &\mapsto y_{\bm{\eta}},
	\end{aligned}
\end{equation}
where $y_{\bm{\eta}}$ is the unique weak solution to \eqref{eq:P-state} associated with $\bm{\eta} = (\eta_1, \eta_2,\ldots, \eta_k)$.

%	\begin{equation}
	%		\label{eq:SE-existence-solution-assumption}
	%		\bm\omega_{\max} \leq 4\pi. 
	%	\end{equation}
%	By \cite{Vazquez1983} (see also \cite[Thm.~4.7]{BresisMarcusPonce2007} and \cite{Bartolucci2005}), for any $\bm{\omega}$ satisfying \eqref{eq:SE-existence-solution-assumption}, the equation \eqref{eq:P-state} admits unique weak solution $y_{\bm{\omega}} \in W^{1,q}_0(\Omega)$ for any $1 \leq q < \frac{N}{N-1} =2$ such that $e^{y_{\bm{\omega}}} \in L^1(\Omega)$. 
%	We shall denote the mapping $(-\infty, 4 \pi]^k \ni \bm{\omega} \mapsto y_{\bm{\omega}} \in W^{1,q}_0(\Omega)$ by $S$. 
The continuity of $S$ and of $\exp(S)$ is demonstrated below.
\begin{lemma}
	\label{lem:control2state-oper-continuity}
	For any $q \in [1,2)$, there exists a  $C = C(q)>0$ satisfying
	\begin{equation}
		\label{eq:SE-apriori-esti}
		\norm{S(\bm{\eta})}_{W^{1,q}_0(\Omega)} \leq C ( \norm{f_0}_{L^p(\Omega)} + \sum_{i =1}^k |\eta_i| )
	\end{equation}
	for any $\bm{\eta} \in \R^k$ satisfying \eqref{eq:SE-existence-solution-assumption-0}. 
	Moreover, one has
	\begin{equation}
		\label{eq:exp-S-esti-L1}
		\int_{\Omega} |e^{S(\bm{\eta})} -1 | dx \leq |\Omega|^{\frac{p-1}{p}} \norm{f_0}_{L^p(\Omega)} + \sum_{i=1}^k |\eta_i|.
	\end{equation}
	Furthermore, there hold
	\begin{align}
		%& \int_{\Omega} [ e^{S(\bm{\omega})} - e^{S(\bm{\eta})} ]^{+} dx \leq \sum_{i=1}^k (\omega_i - \eta_i)^{+}, \quad    
		& \int_{\Omega} |e^{S(\bm{\omega})} - e^{S(\bm{\eta})} | dx \leq \sum_{i=1}^k |\omega_i - \eta_i|, \label{eq:exp-S-continuity-L1-positive-part-whole}
		%& \int_{\Omega} \left|\Delta{S(\bm{\omega})} - \Delta{S(\bm{\eta})} \right| dx \leq 2 \sum_{i=1}^k |\omega_i - \eta_i|, \label{eq:exp-S-Delta-continuity-L1}
		\intertext{and}
		& \norm{S(\bm{\omega}) - S(\bm{\eta})}_{W^{1,q}_0(\Omega)} \leq C \sum_{i =1}^k |\omega_i - \eta_i| \label{eq:S-continuity-W1q}
	\end{align}
%	and
%	\begin{equation}
%		\norm{S(\bm{\omega}) - S(\bm{\eta})}_{W^{1,q}_0(\Omega)} \leq C \sum_{i =1}^k |\omega_i - \eta_i| \label{eq:S-continuity-W1q}
%	\end{equation}
	for any $\bm{\omega}, \bm{\eta} \in \R^k$ with $\bm{\omega}_{\max}, \bm{\eta}_{\max} \leq 4\pi$.
\end{lemma}	
{\it Proof.}  
	The estimate \eqref{eq:exp-S-esti-L1} follows directly from \cite[Cor.~4.B.1]{BresisMarcusPonce2007} and the fact that 
	\begin{align*}
		\norm{f_0 + \sum_{i=1}^k \eta_i \delta_{x_i}}_{\mathcal{M}(\Omega)} & \leq \norm{f_0}_{\mathcal{M}{(\Omega})} + \norm{\sum_{i=1}^k \eta_i \delta_{x_i}}_{\mathcal{M}(\Omega)} = \norm{f_0}_{L^1(\Omega)} + \sum_{i=1}^k |\eta_i| \\
		& \leq |\Omega|^{\frac{p-1}{p}} \norm{f_0}_{L^p(\Omega)} + \sum_{i=1}^k |\eta_i|,
	\end{align*}
	where we have just used the H\"{o}lder inequality to get the last estimate.
	We now rewrite the equation for $S(\bm{\eta})$ as
	\[
	-\Delta S(\bm{\eta}) = f_0 + \sum_{i=1}^k \eta_i \delta_{x_i}  - ( e^{S(\bm{\eta}) } -1 ) \quad \text{in } \Omega, \quad S(\bm{\eta}) = 0 \quad \text{
		on } \Gamma
	\]
	and then apply \cite[Thm.~4.B.1]{BresisMarcusPonce2007} as well as \eqref{eq:exp-S-esti-L1} to have \eqref{eq:SE-apriori-esti}.
	
	Moreover, the  estimate in \eqref{eq:exp-S-continuity-L1-positive-part-whole} follows from  (4.B.20) in \cite{BresisMarcusPonce2007}, upon applying it to \eqref{eq:P-state} with $g(u) := e^u -1$.
	
	Finally, for \eqref{eq:S-continuity-W1q}, by subtracting the equations for $S(\bm{\omega})$ and $S(\bm{\eta})$, we have
	\begin{equation*}
		\left\{
		\begin{aligned}
			-\Delta [ S(\bm{\omega})  - S(\bm{\eta}) ]  & =  \sum_{i=1}^k (\omega_i - \eta_i) \delta_{x_i} - [ e^{ S(\bm{\omega}) }-  e^{ S(\bm{\eta}) } ]&& \text{in } \Omega \\
			S(\bm{\omega})  - S(\bm{\eta}) &= 0 && \text{on } \Gamma.
		\end{aligned}
		\right.
	\end{equation*}
	Applying \cite[Thm.~4.B.1]{BresisMarcusPonce2007} to the above equation, and subsequently employing the estimate in \eqref{eq:exp-S-continuity-L1-positive-part-whole}, we obtain \eqref{eq:S-continuity-W1q}.
\qed   

The following lemma is another version of \Cref{cor:omega-max-4pi} for solutions of the state equation \eqref{eq:P-state}, which will play a role in proving the existence of solutions to the linearization  of the state equation; see \Cref{lem:linearization-state-equation}.
\begin{lemma} 
	\label{lem:omega-max-4pi}
	Assume that $\bm\eta  \in \R^k$ with $\bm{\eta}_{\max} < 4\pi$. Then, there exists a constant $C>0$ independent of $\bm{ \eta}$ satisfying following statements:
	\begin{enumerate}[label=(\alph*)]
		\item 
		If $\bm\eta_{\max} \leq 0$, then there holds
		\begin{equation}
			\label{eq:exponential-state-belong-Linfty-C2S}
			\norm{e^{S(\bm\eta)}}_{L^\infty(\Omega)} \leq e^{C\norm{f_0}_{L^p(\Omega)}}.
		\end{equation}
		\item 
		If $\bm\eta_{\max} \in (0,4\pi)$, then there holds
		\begin{equation}
			\label{eq:exponential-state-belong-L1plus-C2S-norm}
			\norm{e^{S(\bm\eta)}}_{L^{1+\tau}(\Omega)} \leq \pi^{\frac{1}{1+\tau}} e^{C\norm{f_0}_{L^p(\Omega)}}(\frac{2R}{r_0})^{2k}[R^2 +  \frac{k(1+\tau)\bm{\eta}_{\max} r_0^2}{4\pi -(1+\tau)\bm{\eta}_{\max} }]^{\frac{1}{1+\tau}}
		\end{equation}
		for all $\tau \in [0, \frac{4\pi}{\bm\eta_{\max}} -1 )$, where $R$ and $r_0$ are the constants defined in \Cref{prop:Poisson-exponential-esti-many-Dirac}.
	\end{enumerate}
%	\begin{equation}
%		\label{eq:exponential-state-belong-L1plus-C2S}
%		e^{S(\bm\eta)} \in 
%		\left\{
%		\begin{aligned}
%			& L^\infty(\Omega), && \text{if } \bm\eta^{+} = \bm{0},\\
%			& L^{1+\tau}(\Omega) \, \text{for any} \, \tau \in \left[0, \frac{4\pi}{\bm\eta^+_{\max}} -1 \right), && \text{otherwise}.  
%		\end{aligned}
%		\right.
%	\end{equation}
\end{lemma}
{\it Proof.}  	
	Let us put $y:= S(\bm{\eta})$. 
	Setting $y_0 := S(\bm{0})$, that is,
	\[
	-\Delta y_0 + \left( e^{y_0} -1  \right)   = f_0  \quad \text{in } \Omega, \quad
	y_0 = 0 \quad \text{on } \Gamma,
	\]
	%	\begin{equation*}
		%		\left\{
		%		\begin{aligned}
			%			-\Delta y_0 + \left( e^{y_0} -1  \right)  & = f_0  && \text{in } \Omega \\
			%			y_0 &= 0 && \text{on } \Gamma,
			%		\end{aligned}
		%		\right.
		%	\end{equation*}
	yields $y_0 \in C(\overline{\Omega})$, in view of the condition that $p>1$ and of, e.g. \cite[Thm.~4.7]{Troltzsch2010}. 
	There then exists a constant $C>0$ such that
	\begin{equation}
		\label{eq:y0-bounded}
		\norm{y_0}_{C(\overline{\Omega})} \leq C\norm{f_0}_{L^p(\Omega)}.
	\end{equation}
	By putting $y_{+} := S(\bm{\eta}^+)$ and employing the comparison principle; see, e.g. \cite[Cor.~4.B.2]{BresisMarcusPonce2007}, there holds
	\begin{equation}
		\label{eq:y-y0-yplus-comparison}
		y(x), y_0(x) \leq y_{+}(x) \quad \text{for a.e. } x \in \Omega.
	\end{equation}
	Consequently, one has
	\begin{equation}
		\label{eq:exp-y0-yplus-comparison}
		e^{y_{+}(x)} - e^{y_0(x)} \geq 0 \quad \text{for a.e. } x \in \Omega.
	\end{equation}
	Subtracting the equations for $y_{+}$ and $y_0$ yields
	\[
		-\Delta (y_{+} - y_0) + \left( e^{y_{+}} -e^{y_0}  \right)   = \sum_{i=1}^k \eta_i^+ \delta_{x_i} \, \text{in } \Omega \quad
		y_{+} - y_0 = 0 \, \text{on } \Gamma.
%		
%	\left\{
%	\begin{aligned}
%		-\Delta (y_{+} - y_0) + \left( e^{y_{+}} -e^{y_0}  \right)  & = \sum_{i=1}^k \eta_i^+ \delta_{x_i}  && \text{in } \Omega \\
%		y_{+} - y_0 &= 0 && \text{on } \Gamma.
%	\end{aligned}
%	\right.
	\]
	Combining this with \eqref{eq:exp-y0-yplus-comparison} and  the comparison principle gives 
	\[
	y_{+}(x) - y_0(x) \leq z(x) \quad \text{for a.e. } x \in \Omega
	\]
	with $z$ uniquely solving
	\[
	-\Delta z   = \sum_{i=1}^k \eta_i^+ \delta_{x_i}  \quad \text{in } \Omega, \quad
	z = 0 \quad \text{on } \Gamma.
	\]
	%	\[
	%		\left\{
	%		\begin{aligned}
		%			-\Delta z  & = \sum_{i=1}^k \eta_i^+ \delta_{x_i}  && \text{in } \Omega \\
		%			z &= 0 && \text{on } \Gamma.
		%		\end{aligned}
	%		\right.
	%	\]
	Combining this with \eqref{eq:y-y0-yplus-comparison} yields
	\begin{equation}
		\label{eq:y-y0-z-esti}
		y(x) \leq y_{+}(x) \leq y_0(x) + z(x) \quad \text{for a.e. } x \in \Omega.
	\end{equation}
	From this, \eqref{eq:y0-bounded}, and the same arguments as in Cases 1 and 2 in the proof of \Cref{cor:omega-max-4pi}, we obtain \eqref{eq:exponential-state-belong-Linfty-C2S} and \eqref{eq:exponential-state-belong-L1plus-C2S-norm}.
\qed   
The number $4\pi$ is critical as it is seen in the following.
\begin{proposition}
	\label{prop:critical-4pi}
	Let  $\bm\eta \in \R^k$ be arbitrary, but fixed. Then, the following assertions are valid.
	\begin{enumerate}[label=(\alph*)]
		\item \label{item:4pi}  If $\bm\eta_{\max} = 4 \pi$, then, for any $\tau >0$, $e^{S(\bm\eta)}$ does not belong to $L^{1 + \tau}(\Omega)$.
		\item \label{item:less-4pi} If $\bm{\eta}_i < 4\pi$ for some $1 \leq i \leq k$, then there exists a constant $\rho >0$ such that the restriction of $e^{S(\bm\eta)}$ to $B_{\R^2}(x_i,\rho)$ belongs to $L^{1+\tau_i}(B_{\R^2}(x_i,\rho))$ for some $\tau_i > 0$.
	\end{enumerate}
\end{proposition}
{\it Proof.}  
	We now set $y:= S(\bm\eta)$.

	\noindent{\it Ad: \ref{item:4pi}}.
	Let $\tau >0$ be arbitrary, but fixed. Assume that, $\eta_{i_0} = 4\pi$ for some $i_0 \in \{1,2,\ldots,k\}$. 
	Setting  $w := e^y -1 - f_0$ yields $w \in L^1(\Omega)$. 
	Moreover, there exists a constant $\epsilon_0 >0$ such that $B_{\R^2}(x_{i}, \epsilon_0)$, $1 \leq i \leq k$, are disjoint open balls in $\Omega$. Consequently, there holds
	\begin{equation}
		\label{eq:state-in-vinicity-xi}
		-\Delta y + w = 4\pi \delta_{x_{i_0}} \quad \text{in} \quad \mathcal{D}'(B_{\R^2}(x_{i_0}, \epsilon_0)).
	\end{equation}	
	We now define the angular mean $\tilde{y}: (0,\epsilon_0) \to \R$  of $y$ as
	\[
	\tilde{y}(r) = \frac{1}{2\pi} \int_0^{2 \pi} y(r, \theta) d\theta,
	\]
	where $(r,\theta)$ denotes the polar coordinates in $\R^2$. Since $w$ has no point mass at $x_{i_0}$, we deduce from \cite[Cor.~2]{Vazquez1983} that
	$
	\lim\limits_{r \to 0^+} \frac{\tilde{y}(r)}{ \ln(1/r)} = \frac{4 \pi}{2 \pi} = 2.
	$
	There then exists a constant $\epsilon \in (0, \epsilon_0)$ such that
	\[
	\tilde{y}(r) \geq (2 - \frac{2\tau}{1+\tau}) \ln (\frac{1}{r}) = \frac{2}{1+\tau} \ln (\frac{1}{r}) \quad \text{for all } r \in (0, \epsilon).
	\]
	We thus deduce from this, together with Jensen's inequality for integrals, that
	\begin{multline*}
		\frac{1}{2\pi}\int_{B_{\R^2}(x_{i_0}, \epsilon)} e^{(1+\tau)y} dx  = \int_{0}^{\epsilon} r dr \left[ \frac{1}{2\pi} \int_{0}^{2\pi} e^{(1+\tau){y}(r, \theta)} d\theta \right] \\
		\begin{aligned}
			& \geq \int_{0}^{\epsilon} r dr \exp\left[ \frac{1}{2\pi} \int_0^{2\pi} {(1+\tau){y}(r, \theta)} d\theta \right] 
			%\\
			= \int_{0}^{\epsilon} r \exp\left[(1+\tau) \tilde{y}(r)\right] dr \\
			& \geq \int_{0}^{\epsilon} r \times \frac{1}{r^2} dr = \infty.
		\end{aligned}
	\end{multline*}
	The derised conclusion of assertion \ref{item:4pi} then follows.
	
	\noindent{\it Ad: \ref{item:less-4pi}}. 
	It suffices to consider the case, where $\eta_i \in (0, 4\pi)$. Let $y_0$ and $\tilde{y}$ be unique solutions to
	\begin{equation*}
		\left\{
		\begin{aligned}
			-\Delta y_0  & = |f_0|   && \text{in } \Omega \\
			y_0 &= 0 && \text{on } \Gamma
		\end{aligned}
		\right. \quad \text{and} \quad 
		\left\{
		\begin{aligned}
			-\Delta \tilde{y}   & =  \sum_{i =1}^k \eta_i^{+} \delta_{x_{i}}  && \text{in } \Omega \\
			\tilde{y} &= 0 && \text{on } \Gamma,
		\end{aligned}
		\right.
	\end{equation*}
	respectively. Since $f_0 \in L^p(\Omega)$ with $p > 1 = \frac{N}{2}$, then $y_0 \in L^\infty(\Omega)$. Moreover, by comparison principle; see, e.g. \cite[Cor.~4.B.2]{BresisMarcusPonce2007}, there hold
	$
	y_0, \tilde{y} \geq 0  
	$
	and 
	$
	y \leq y_0 + \tilde{y}
	$
	a.e. in $\Omega$; see \eqref{eq:y-y0-z-esti}. Consequently, we have for a.e. $x \in \Omega$ that
	\begin{equation}
		\label{eq:y-exp-bounded-above}
		e^{y(x)} \leq Ce^{\tilde{y}(x)} \quad \text{for some constant } C>0.
	\end{equation}
	On the other hand, by \eqref{eq:exponential-esti-Dirac-many-Bxj}, one has 
	$
	\int_{B_{\R^2}(x_i,\rho)} \exp\left[\frac{(4\pi - \alpha) \tilde{y}}{\eta_i}\right]dx  \leq C_\alpha
	$
	for some constant $\rho$ and for all $\alpha \in (0, 4\pi)$.
	From this and \eqref{eq:y-exp-bounded-above}, we have the desired conclusion. 
\qed

%%% The continuity of Nemytskii operator
The continuity on $L^{1+\tau}$-space of Nemytskii's operator induced by the exponential function is shown in the next lemma.
\begin{lemma}
	\label{eq:exponential-Nemytskii-continuity}
	Let $\bm\omega  \in \R^k$ and 
	$\tau >0$ satisfy $(1+ \tau ) \bm\omega_{\max} < 4 \pi$.
%	\begin{equation}
%		\label{eq:tau-4pi-condition}
%		(1+ \tau ) \bm\omega_{\max} < 4 \pi.
%	\end{equation}
	Then,
	there holds 
	\begin{equation}
		\label{eq:exp-Nemytskii-continuity}
		e^{S(\bm\omega + \bm\eta)} \to e^{S(\bm\omega)} \quad \text{in} \quad L^{1+\tau}(\Omega) \quad \text{as} \quad \|\bm\eta\|_{\infty} \to 0.
	\end{equation}
\end{lemma}
{\it Proof.}  
	%Assume that \eqref{eq:tau-4pi-condition} is fulfilled.
	%By \Cref{lem:omega-max-4pi}, there holds $e^{S(\bm\omega + \bm\eta)}, e^{S(\bm\omega)} \in L^{1 +\tau}(\Omega)$, provided $\|\bm\eta\|_{\infty}$ is small enough.
	Let $\tau_0 > \tau$ be such that $(1+ \tau_0 ) \bm\omega_{\max} < 4 \pi$.
	We then deduce from \eqref{eq:exponential-state-belong-Linfty-C2S} and \eqref{eq:exponential-state-belong-L1plus-C2S-norm} that there exists a constant $\epsilon_0  >0$ such that the family
	\[
		\{ e^{S(\bm\omega + \bm\eta)} \, | \, \bm{ \eta} \in \R^k,  \|\bm\eta\|_{\infty} < \epsilon_0 \}
	\]
	is uniformly bounded in $L^{1 +\tau_0}(\Omega)$.
	On the other hand, we conclude from the continuity of the control-to-state operator $S$ (see  \Cref{{lem:control2state-oper-continuity}}) that
	$S(\bm\omega + \bm\eta) \to S(\bm\omega)$ in measure in $\Omega$ as $\|\bm\eta\|_{\infty} \to 0$.
%	\begin{equation*}
%		%\label{eq:S-eta-limit}
%		S(\bm\omega + \bm\eta) \to S(\bm\omega) \quad \text{in measure in } \Omega \quad \text{as} \quad \|\bm\eta\|_{\infty} \to 0.
%	\end{equation*}
	We then apply \cite[Chap.~XI, Prop.~3.10]{Visintin1996} to derive \eqref{eq:exp-Nemytskii-continuity}.
\qed
The following theorem on the continuous differentiability of the control-to-state mapping $S$ on the open set $\mathcal{O}_{\bm{b}}$ is the main result of this section.
\begin{theorem}
	\label{thm:F-diff-S}
	Assume that $0 < {\bm{b}}_{\max} < 4 \pi$. Then, the control-to-state operator $S$ is of class $C^1$ on $\mathcal{O}_{\bm{b}}$. Moreover, for any $\bm\eta \in \mathcal{O}_{\bm{ b } }$ and $\bm\omega\in \R^k$, the derivative $z := S'(\bm\eta)\bm\omega$ is characterized via
	\begin{equation}
		\label{eq:derivative-S-1st}
		-\Delta z + e^{S(\bm\eta)} z   = \sum_{i=1}^{k} \omega_i \delta_{x_i} \, \text{in } \Omega, \quad
		z = 0 \, \text{on } \Gamma.
%		
%		\left\{
%		\begin{aligned}
%			-\Delta z + e^{S(\bm\eta)} z  & = \sum_{i=1}^{k} \omega_i \delta_{x_i} && \text{in } \Omega, \\
%			z &= 0 && \text{on } \Gamma.
%		\end{aligned} 
%		\right.
	\end{equation}
\end{theorem}
{\it Proof.}  
	Let us introduce the mapping defined on an open subset of $W_{r_{\bm{b}}} \times \R^k$ as follows
	\[
	\begin{aligned}
		F:  V_{\bm{b}} \times \mathcal{O}_{\bm{b}} & \to \mathcal{D} + L^{r_{\bm{b}}}(\Omega)  \\
		(y,\bm \eta) &\mapsto -\Delta y + e^y - 1 - f_0 - \sum_{i=1}^{k} \eta_i \delta_{x_i},
	\end{aligned}
	\]
	where the open sets $\mathcal{O}_{\bm{b}}$ and $V_{\bm{b}}$ are defined in \eqref{eq:Ob-set} and \eqref{eq:Vb-set}, respectively, while the space $\mathcal{D}$ is given in \eqref{eq:D-space}.
	Obviously, there holds
	\[
	F(y, \bm{\eta}) = F_0(y) + H(\bm{\eta}) \quad \text{for all } (y,\bm \eta) \in V_{\bm{b}} \times \mathcal{O}_{\bm{b}}
	\]
	with $F_0$ being defined as in \eqref{eq:F0-operator} and $H: \mathcal{O}_{\bm{b}}  \to \mathcal{D} + L^{r_{\bm{b}}}(\Omega)$ being given by
	\[
	H(\bm{\eta}) := - f_0 - \sum_{i=1}^{k} \eta_i \delta_{x_i} \quad \text{for } \bm{\eta} \in \mathcal{O}_{\bm{b}}.
	\]
	%	
	%	
	%	By definition, for any $y \in V_{\bm{b}}$, there exist $g \in L^{r_{\bm{b}}}(\Omega)$ and $\bm{\eta} \in \R^k$ with $\eta_i < b_i$ for all $1 \leq i \leq k$ such that $y \in W^{1,1}_0(\Omega)$ and
	%	\[
	%		-\Delta y = g + \sum_{i=1}^k \eta_i \delta_{x_i} \quad \text{in } \Omega.
	%	\] 
	%	By \Cref{cor:omega-max-4pi}, one has $e^y \in L^{1+ \tau}(\Omega)$ for any $\tau \in \left[0, \frac{4 \pi}{\bm{\eta}_{\max}} -1\right)$. Since $\bm{\eta}_{\max} < \bm{b}_{\max}$, there holds
	%	$e^y \in L^{\frac{4 \pi}{\bm{b}_{\max}}}(\Omega)$ and thus $e^y - 1 - f_0 \in  L^{\min\{ r_{\bm{b}}, p\}}(\Omega)$. This shows the well-definedness of $F$. 
	%	
	By \Cref{lem:F0-differentiability}, $F(\cdot, \bm\eta)$ is of class $C^1$ over $V_{\bm{b}}$. Moreover, in view of  \Cref{lem:F0-der-isomorphism}, $\frac{\partial F}{\partial y}(y, \bm{\eta}) = \frac{\partial F_0}{\partial y}(y)$ is an isomorphism as a mapping from $W_{r_{\bm{b}}}$ to $\mathcal{D} + L^{r_{\bm{b}}}(\Omega)$. 
	As an immediate consequence of the implicit function theorem; see, e.g. \cite[Thm.~4.B]{Zeidler_vol1}, we can conclude that $F$ is of class $C^1$.
	Moreover,
	a simple computation, together with \eqref{eq:F0-1st-der}, indicates that
	$S'(\bm\eta)$ is determined by \eqref{eq:derivative-S-1st}.
	%	 Moreover,
	%	a simple computation, together with \eqref{eq:F0-1st-der} and \eqref{eq:F0-2nd-der}, indicates that
	%	$S'(\bm\eta)$ and $S''(\bm{\eta})$ are determined by \eqref{eq:derivative-S-1st} and \eqref{eq:derivative-S-2nd}, respectively.
\qed   

%%% The remark on the critical case omega_max = 4 \pi
\begin{remark}
	\label{rem:critical-case-4pi}	
	If $\bm{\eta}_{\max} = 4\pi$, then the differentiability of the control-to-state operator $S$ may fail to hold. The reason for this is that $e^{S(\bm{\eta})}$ belongs to $L^1(\Omega)$ only (see assertion \ref{item:4pi} in \Cref{prop:critical-4pi}). As a result, \Cref{lem:linearization-state-equation} cannot be applied, and thus the linearized state equation \eqref{eq:linearization-state-equation} may not admit solutions in $W^{1,1}_0(\Omega)$.  
	The existence of weak solutions to \eqref{eq:linearization-state-equation} in the case where $e^{S(\bm{\eta})} \in L^1(\Omega)$, with the right-hand side belonging to $L^s(\Omega)$ with $s > 1$, was considered in  \cite[Thms.~1.1 and ~3.3]{Lee2025}.
	%\Cref{prop:elliptic-L1-zero-order-term} below.
\end{remark}

%%%% Necessary optimality conditions
\section{Reduced objective functional and proofs of main results}
\label{sec:proof-mainresult}

\subsection{Reduced objective functional}
\label{subsec:reduced-cost-func}
%From \Cref{rem:restriction-feasible-set}, from now on, we can assume that $\bm{0} \leq \bm{M} \leq \bm{4\pi}$, that is,
%\[
%{0} \leq {M}_i \leq {4\pi} \quad \text{for all } 1 \leq i \leq k.
%\]
By \eqref{eq:control-2-state-oper}, the control-to-state mapping $S$ is well-defined on $(-\infty, 4\pi]^k \subset \R^k$. Consequently, the objective functional $J$ of \eqref{eq:P-original} can be reduced as
\begin{equation}
	\label{eq:objective-func-reduced}
	\hat J(\bm{\eta}) := J(S(\bm{\eta}), \bm{\eta}) = \frac{1}{2} \norm{S(\bm{\eta}) - y_d}^2_{L^2(\Omega)} + \kappa \sum_{i =1}^k |\eta_i|, \quad \bm{ \eta} \in  (-\infty, 4\pi]^k.
\end{equation}
%%%%% The objective functional J(\eta) = T(\eta) + \kappa j(\eta) with j(\eta) = \| \eta \|_{l^1}
Let us decompose the reduced objective functional $\hat J$ into two summands as 
\[
\hat J(\bm{\eta}) = T(\bm{\eta}) + \kappa j(\bm{\eta}),
\]
where
\begin{equation} \label{eq:T-func-j-func}
	T(\bm{\eta})  = \frac{1}{2} \norm{S(\bm{\eta}) - y_d}^2_{L^2(\Omega)}   \quad
	\text{and} \quad
	j(\bm{\eta}) = \sum_{i =1}^k |\eta_i|. %%%\sum_{i =1}^k j_0(\eta_i)  \label{eq:j-func} \quad \text{with } j_0(t) = |t| \, \text{for} \quad t \in \R.
\end{equation}
Obviously, $j$ is convex and globally Lipschitz continuous.

As already  shown in \Cref{thm:F-diff-S}, the control-to-state mapping $S$ is of class $C^1$ on the open set $\mathcal{O}_{\bm{b}} \subset (-\infty, 4 \pi]^k$ with $\bm{b}_{\max} \in (0, 4\pi)$.
Therefore, by the chain rule and by the definition of $T$ in \eqref{eq:T-func-j-func}, the functional $T$ is also of class $C^1$ on this open set, as stated below.  
\begin{theorem}
	\label{thm:T-classC2}
	For any $\bm{b} \in \R^k$ with $\bm{b}_{\max} \in (0, 4\pi)$, the functional $T$, defined in \eqref{eq:T-func-j-func}, is of class $C^1$ over the open set $\mathcal{O}_{\bm{b}}$, defined as in \eqref{eq:Ob-set}. Moreover, for any $\bm{\eta} \in \mathcal{O}_{\bm{b}}$ and $\bm{\omega} \in \R^k$, there holds
	\begin{align}
		T'(\bm{\eta})(\bm{\omega}) & = \int_{\Omega} (S(\bm{\eta}) - y_d)z_{\bm{\eta}; \bm{\omega}} dx  = \sum_{i=1}^k \varphi_{\bm{\eta}}(x_i) \omega_i, \label{eq:T-der-1st}		
	\end{align} 
	where $z_{\bm{\eta}; \bm{\omega}} := S'(\bm{\eta})\bm{\omega}$ and $\varphi_{\bm{\eta}}$ is the unique solution in $W^{1,s_*}_0(\Omega)$ 
	for some $s_* >2$
	%	for all $s \in \R$ satisfying
	%	\begin{equation}
		%		\label{eq:s-exponent-condition}
		%		\frac{1}{2} > \frac{1}{s} > \frac{\bm{b}_{\max}}{4\pi} - \frac{1}{2}
		%	\end{equation}
	to 
	\begin{equation}
		\label{eq:adjoint-state}
		-\Delta \varphi_{\bm{\eta}} + e^{S(\bm{\eta})} \varphi_{\bm{\eta}}  = S(\bm{\eta}) - y_d \, \text{in } \Omega \quad
		\varphi_{\bm{\eta}}  = 0 \, \text{on } \Gamma.
%		
%		\left\{
%		\begin{aligned}
%			-\Delta \varphi_{\bm{\eta}} + e^{S(\bm{\eta})} \varphi_{\bm{\eta}}& = S(\bm{\eta}) - y_d && \text{in } \Omega \\
%			\varphi_{\bm{\eta}} &= 0 && \text{on } \Gamma.
%		\end{aligned}
%		\right.
	\end{equation}
	%	Furthermore, for any $\bm{\eta} \in \mathcal{O}_{\bm{b}}$ and $\bm{\omega}, \bm{\psi} \in \R^k$, there holds
	%	\begin{align}
		%		T''(\bm{\eta})(\bm{\omega}, \bm{\psi}) & = \int_{\Omega} (S(\bm{\eta}) - y_d)w_{\bm{\eta}; \bm{\omega}, \bm{\psi}} + z_{\bm{\eta}; \bm{\omega}}z_{\bm{\eta}; \bm{\psi}} dx \notag \\
		%		&  = \int_{\Omega} \left[1 - e^{S(\bm{\eta})} \varphi_{\bm{\eta}} \right] z_{\bm{\eta}; \bm{\omega}} z_{\bm{\eta}; \bm{\psi}} dx \label{eq:T-der-2nd}		
		%	\end{align}
	%	with $w_{\bm{\eta}; \bm{\omega}, \bm{\psi}} := S''(\bm{\eta})(\bm{\omega}, \bm{\psi})$.
\end{theorem}
{\it Proof.}  
	In view of \Cref{thm:F-diff-S} and the chain rule, $T$ is of class $C^1$ over $\mathcal{O}_{\bm{b}}$. Moreover, we obviously derive the first identity in \eqref{eq:T-der-1st}. 
	% and \eqref{eq:T-der-2nd}. 
	For the second one, from the fact that $\bm{b}_{\max} \in (0, 4\pi)$, we can choose a number $s_* \in (2, 4]$ satisfying
	\begin{equation}
		\label{eq:s-exponent-condition}
		\frac{1}{s_*} > \frac{\bm{b}_{\max}}{4\pi} - \frac{1}{2}.
	\end{equation}
	We now prove the $W^{1,s_*}_0(\Omega)$-regularity 
	of solutions to \eqref{eq:adjoint-state}. 		
	To this end, 
	we  employ \Cref{lem:omega-max-4pi} to arrive at 
	\begin{equation}
		\label{eq:state-belong-Lr}
		e^{S(\bm{\eta})} \in L^{\frac{4 \pi}{\bm{b}_{\max}}}(\Omega) \quad \text{for all} \quad \bm{\eta} \in \mathcal{O}_{\bm{b}}.
	\end{equation}
	Let us take $\tau \in \R$ such that
	$
	\frac{1}{\tau} + \frac{\bm{b}_{\max}}{4 \pi} = \frac{1}{s_*} + \frac{1}{2}.
	$
	Since $s_*>2$, there holds $\tau > 1$. Moreover, due to \Cref{lem:linearization-state-equation} and \eqref{eq:state-belong-Lr}, as well as the fact that $\frac{4 \pi}{\bm{b}_{\max}} >1$,  \eqref{eq:adjoint-state} admits a unique solution $\varphi_{\bm{\eta}}$, which belongs to $W^{1,q}_0(\Omega)$ for all $q \in [1,2)$. The continuous embedding $W^{1,q}_0(\Omega) \hookrightarrow L^{\tau}(\Omega)$ for some $q$ sufficiently close to $2$ then implies that $\varphi_{\bm{\eta}} \in L^{\tau}(\Omega)$. Combining this with \eqref{eq:state-belong-Lr}, the H\"{o}lder inequality, and the choice of $\tau$ yields
	\[
	e^{S(\bm{\eta})} \varphi_{\bm{\eta}} \in L^{\frac{2s_*}{s_*+2}}(\Omega),
	\]
	which, along with the fact that $S(\bm{\eta}) \in L^{r}(\Omega)$ for all $r \geq 1$, according to the embedding $W^{1,q}_0(\Omega) \hookrightarrow L^{r}(\Omega)$ for some $q$ sufficiently close to $2$, then gives
	\[
	S(\bm{\eta}) - y_d - e^{S(\bm{\eta})} \varphi_{\bm{\eta}} \in L^{\frac{2s_*}{s_*+2}}(\Omega) \hookrightarrow W^{-1,s_*}(\Omega).
	\]
	Here we have just used the Sobolev embedding theorem to derive the last embedding. Applying \cite[Thm.~0.5]{JerisonKenig1995} to \eqref{eq:adjoint-state}
	%	\begin{equation*}
		%		\left\{
		%		\begin{aligned}
			%			-\Delta \varphi_{\bm{\eta}}  & = S(\bm{\eta}) - y_d - e^{S(\bm{\eta})} \varphi_{\bm{\eta}} && \text{in } \Omega \\
			%			\varphi_{\bm{\eta}} &= 0 && \text{on } \Gamma
			%		\end{aligned}
		%		\right.
		%	\end{equation*}
	yields $\varphi_{\bm{\eta}} \in W^{1,s_*}_0(\Omega)$.
	
	We now show the second identity in \eqref{eq:T-der-1st}. 
	% and \eqref{eq:T-der-2nd}. 
	To achieve this, we  test \eqref{eq:derivative-S-1st} and \eqref{eq:adjoint-state} by $\varphi_{\bm{\eta}}$ and  $z_{\bm{\eta}; \bm{\omega}}$, respectively, to have
	\[
		\left\{
		\begin{aligned}
			& \int_\Omega \nabla \varphi_{\bm{\eta}} \cdot \nabla z_{\bm{\eta}; \bm{\omega}} + e^{S(\bm{\eta})} \varphi_{\bm{\eta}} z_{\bm{\eta}; \bm{\omega}} dx = \sum_{i=1}^k \varphi_{\bm{\eta}}(x_i) \omega_i \\
			%\intertext{and}
			& \int_\Omega \nabla \varphi_{\bm{\eta}} \cdot \nabla z_{\bm{\eta}; \bm{\omega}} + e^{S(\bm{\eta})} \varphi_{\bm{\eta}} z_{\bm{\eta}; \bm{\omega}} dx = \int_{\Omega} (S(\bm{\eta}) - y_d)z_{\bm{\eta}; \bm{\omega}} dx.
		\end{aligned}
		\right.
	\]
	Subtracting these equations then yields the second identity in \eqref{eq:T-der-1st}.
	%	Similarly, testing \eqref{eq:derivative-S-2nd} and \eqref{eq:adjoint-state} by $\varphi_{\bm{\eta}}$ and  $w_{\bm{\eta}; \bm{\omega}, \bm{\psi}}$, respectively, and then subtracting the obtained equations thus gives the second identity in \eqref{eq:T-der-2nd}.
\qed   

\subsection{Proof of \Cref{thm:1st-OCs-P-orig} under assumption \ref{ass:regularcase}}
\label{subsec:proof-main-regular}

%\noindent{\bf Proof of \Cref{thm:1st-OCs-P-orig} under assumption \ref{ass:regularcase}.}
\noindent{\it Proof.}
	Assume that \ref{ass:regularcase} is fulfilled. We then obtain from \Cref{thm:minimizer-existence} that $\bar{\bm{\eta}}_{\max} < 4\pi$.
	The functional $T$, defined in \eqref{eq:T-func-j-func}, is of class $C^1$ in a neighborhood of $\bm{\bar \eta}$, by virtue of \Cref{thm:T-classC2}. 
	Obviously, $\bar y$ and $\bar{\bm{\eta}}$ satisfy \eqref{eq:OC-1st-state-P-orig}. Let $\bar{\varphi}$ be the unique solution in $W^{1, \bar s}_0(\Omega) \hookrightarrow C(\overline{\Omega})$ for some $\bar s >2$ to \eqref{eq:OC-1st-adjoint-state-P-orig}, which is well-defined as in \Cref{thm:T-classC2}.	
In order to show \eqref{eq:OC-1st-variational-P-orig}, we first take any $\bm{\eta} = (\eta_1, \eta_2, \ldots, \eta_k) \in \R^k$.
There thus exists a constant $\rho_0 \in (0,1)$ satisfying 
\[
\bar \eta_j + \rho(\eta_j - \bar \eta_j) \leq 4\pi \quad \text{for all } 1 \leq j \leq k \quad \text{and for all } 0 < \rho < \rho_0. 
\]
Therefore, $\bar{ \bm{\eta}} + \rho (\bm{\eta} - \bar{ \bm{\eta}})$, with $0 < \rho < \rho_0$, is an admissible control for \eqref{eq:P-original} and we thus have
\[
\hat J(\bar{ \bm{\eta}} + \rho (\bm{\eta} - \bar{ \bm{\eta}})) \geq \hat J(\bar{\bm{\eta}}),
\] 
by the optimality of $\bar{\bm{\eta}}$. 
From this and the convexity of $j$, there holds
\begin{align*}
	0 & \leq \frac{1}{\rho} \left[  \hat J(\bar{ \bm{\eta}}     + \rho (\bm{\eta} - \bar{ \bm{\eta}}    )) - \hat J(\bar{\bm{\eta}}    ) \right] \\
	& = \frac{1}{\rho} \left[  T(\bar{ \bm{\eta}}     + \rho (\bm{\eta} - \bar{ \bm{\eta}}    )) - T(\bar{\bm{\eta}}    ) \right] + \frac{\kappa}{\rho} \left[  j(\bar{ \bm{\eta}}     + \rho (\bm{\eta} - \bar{ \bm{\eta}}    )) - j(\bar{\bm{\eta}}    ) \right] \\
	& \leq \frac{1}{\rho} \left[  T(\bar{ \bm{\eta}}     + \rho (\bm{\eta} - \bar{ \bm{\eta}}    )) - T(\bar{\bm{\eta}}    ) \right] + \kappa \left[ j(\bm{\eta}) - j(\bm{\bar \eta}    )\right]
\end{align*}
for all $0< \rho < \rho_0$. By letting now $\rho \to 0^+$ and using \Cref{thm:T-classC2}, we have
\begin{equation}
	\label{eq:1st-OC-direction-active}
	0  \leq T'(\bm{\bar \eta}    )(\bm{\eta} - \bm{\bar \eta}    ) + \kappa \left[ j(\bm{\eta}) - j(\bm{\bar \eta}    )\right] = \sum_{i=1}^{k}  \bar{\varphi}    (x_i)(\eta_i - \bar{\eta}_i    ) + \kappa \left[ |{\eta}_i| - |{\bar \eta}_i    | \right]
\end{equation}
for all $\bm{\eta} \in \R^k$. A straightforward argument then implies \eqref{eq:OC-1st-variational-P-orig}.
\qed

\subsection{Proof of \Cref{thm:1st-OCs-P-orig} under assumption \ref{ass:irrergularcase}}
\label{subsec:proof-main-irregular}
	
	Under assumption \ref{ass:irrergularcase}, we have $e^{\bar y} \in L^1(\Omega)$ only (see \Cref{prop:critical-4pi}), and the control-to-state operator $S$ may fail to be differentiable at $\bm{\bar \eta}$ (see \Cref{rem:critical-case-4pi}). Therefore, the argument establishing the main results is more involved than the one in \Cref{subsec:proof-main-regular}.
	
\medskip

%\noindent{\bf Proof of \Cref{thm:1st-OCs-P-orig} under assumption \ref{ass:irrergularcase}.}	
\noindent{\it Proof.}	
	Assume that \ref{ass:irrergularcase} is satisfied.
	Clearly, \eqref{eq:OC-1st-state-P-orig} is verified. 
	Let $\bar{\varphi} \in H^1_0(\Omega) \cap L^\infty(\Omega)$ denote the unique solution to \eqref{eq:OC-1st-adjoint-state-P-orig}, which is determined by \cite[Thms.~1.1 and~3.3]{Lee2025} in the case $e^{\bar y} \in L^1(\Omega)$.

	It remains to show \eqref{eq:OC-1st-variational-P-orig} and the local continuity of $\bar\varphi$. To this end, we first consider a family of regularized optimal control problems \eqref{eq:P-regularized} associated with \eqref{eq:P-original}, then derive the corresponding optimality systems, and finally pass to the limit as $\epsilon \to 0^+$. These steps are carried out in the three parts below.
	
	\medskip 
	
	\noindent{\it Step 1: Regularized optimal control problems.}
	Since $(\bar y,\bar{\bm{\eta}})$ is a strict local minimizer, it follows from \Cref{rem:restriction-feasible-set} that an $\epsilon_0 >0$ exists and satisfies 
	\begin{equation}
		\label{eq:strict-minimizer}
		\hat J(\bm{\bar \eta}) < \hat{J}(\bm{\eta}) \quad \forall \bm{ \eta} \in \R^k \backslash\{\bm{\bar \eta} \} \, : \bm{ \eta}_{\max} \leq 4\pi,  |\eta_i - \bar \eta_i| \leq \epsilon_0 \, \forall 1 \leq i \leq k.
	\end{equation}
	Let us set $\bm{a} := \bm{\bar \eta} - \epsilon_0 \bm{1}$. 
	For any $\epsilon \in (0,\epsilon_0)$, we define a vector $\bm{b}^\epsilon \in \R^k$, depending on $\bm{\bar \eta}$, as follows
	\begin{equation}
		\label{eq:ab-def}
		\left\{
		\begin{aligned}
			&b_i^\epsilon = b_i^{\epsilon_0} := \bar{\eta}_i + \epsilon_0 && \text{if } i \notin I_{4\pi}(\bm{\bar \eta}),\\
			&b_i^\epsilon := 4\pi - \epsilon  && \text{if } i \in I_{4\pi}(\bm{\bar \eta}).
		\end{aligned}
%		\begin{aligned}
%			& a_i := \bar{\eta}_i - \epsilon_0, \quad b_i^\epsilon = b_i^{\epsilon_0} := \bar{\eta}_i + \epsilon_0 && \text{if } i \notin I_{4\pi}(\bm{\bar \eta}),\\
%			& a_i := 4\pi - \epsilon_0, \quad b_i^\epsilon := 4\pi - \epsilon && \text{if } i \in I_{4\pi}(\bm{\bar \eta}).
%		\end{aligned}
		\right.
	\end{equation}
	We can reduce $\epsilon_0$ if necessary to have for all $\epsilon \in (0, \epsilon_0)$ and $1 \leq i \leq k$ that
	\begin{equation}
		\label{eq:ab-epsilon-zero}
		\bm{b}^\epsilon_{\max} < 4\pi \quad \text{and} \quad
		\left\{
		\begin{aligned}
			& b_i^\epsilon > a_i > 0 && \text{if } \bar \eta_i >0,\\
			& 0 > b_i^\epsilon > a_i && \text{if } \bar \eta_i < 0.
		\end{aligned}
		\right.
	\end{equation}
Obviously, one has $\bm{a}_{\max}, \bm{b}^\epsilon_{\max} < 4\pi$ for all $\epsilon \in (0, \epsilon_0)$. 
We then consider regularized optimal control problems, depending on $\epsilon$, of \eqref{eq:P-original} as follows
\begin{equation}
	\label{eq:P-regularized}
	%\tag{\text{$P_{\bm{M}}$} }
	\tag{\text{${\rm P_{\epsilon}}$}}
	\min \{  \hat J(\bm{\eta}) =J(S(\bm{ \eta}), \bm{\eta}) \, | \,  \bm{\eta} \in  \R^k,  \, \bm{a} \leq \bm{ \eta} \leq \bm{b}^\epsilon \}.
\end{equation}
By a standard argument; see the proof of \Cref{thm:minimizer-existence},  \eqref{eq:P-regularized} admits at least one global minimizer $\bar{\bm{\eta}}^\epsilon$.

\medskip
\noindent{\it Step 2: Deriving the optimality conditions for \eqref{eq:P-regularized}.}
	We now show that there exists 
	there exists an adjoint state $\bar{\varphi}^{\epsilon} \in W^{1,s_\epsilon}_0(\Omega) \hookrightarrow H^1_0(\Omega) \cap L^\infty(\Omega)$ for some $s_\epsilon>2$, that, together with  $\bar{\bm{\eta}}^{\epsilon}$ and $\bar y^\epsilon := S(\bar{\bm{\eta}}^\epsilon)$, satisfies the system
	\begin{subequations}
		\label{eq:OCs-1st}
		\begin{align}
			&
			\label{eq:OC-1st-state} 
			-\Delta \bar y^{\epsilon} + (e^{\bar y^{\epsilon}} - 1) = f_0(x) + \sum_{i =1}^k \bar\eta_i^{\epsilon} \delta_{x_i} \,  \text{in } \Omega, \quad
			\bar y^{\epsilon} = 0 \, \text{on } \Gamma, \\
%			\left\{
%			\begin{aligned}
%				-\Delta \bar y^{\epsilon} + (e^{\bar y^{\epsilon}} - 1) &= f_0(x) + \sum_{i =1}^k \bar\eta_i^{\epsilon} \delta_{x_i} && \text{in } \Omega, \\
%				\bar y^{\epsilon} &= 0 && \text{on } \Gamma,
%			\end{aligned}
%			\right.  \label{eq:OC-1st-state} \\
			& 
			 \label{eq:OC-1st-adjoint-state} 
			-\Delta \bar \varphi^{\epsilon} + e^{\bar y^{\epsilon}} \bar\varphi^{\epsilon} = \bar y^{\epsilon} - y_d \, \text{in } \Omega, \quad
			\bar \varphi^{\epsilon} = 0 \, \text{on } \Gamma,
%			
%			\left\{
%			\begin{aligned}
%				-\Delta \bar \varphi^{\epsilon} + e^{\bar y^{\epsilon}} \bar\varphi^{\epsilon} &= \bar y^{\epsilon} - y_d && \text{in } \Omega, \\
%				\bar \varphi^{\epsilon} &= 0 && \text{on } \Gamma,
%			\end{aligned}
%			\right. \label{eq:OC-1st-adjoint-state} 
			\intertext{and}
			& 
			\left\{
			\begin{aligned}
				& - \frac{\bar\varphi^{\epsilon}(x_i)}{\kappa} \geq \sign(\bar\eta_i^{\epsilon}) && \text{if } \bar\eta_i^{\epsilon} = b_i^\epsilon, \\
				& - \frac{\bar\varphi^{\epsilon}(x_i)}{\kappa} \leq \sign(\bar\eta_i^{\epsilon}) && \text{if } \bar\eta_i^{\epsilon} = a_i, \\
				& - \frac{\bar\varphi^{\epsilon}(x_i)}{\kappa} = \sign(\bar\eta_i^{\epsilon}) && \text{if } \bar\eta_i^{\epsilon} \in (a_i, b_i^\epsilon)\backslash\{0\}, \\ 
				& - \frac{\bar\varphi^{\epsilon}(x_i)}{\kappa} \in [-1,1] && \text{if } \bar\eta_i^{\epsilon} = 0, a_i < 0 < b_i^\epsilon, \\
			\end{aligned}
			\right.		
			%			\quad	
			%			\left\{
			%			\begin{aligned}
				%				& - \frac{\bar\varphi^{\epsilon}(x_i)}{\kappa} \geq 1 && \text{if } \bar\eta_i^{\epsilon} = b_i^\epsilon > 0, \\
				%				& - \frac{\bar\varphi^{\epsilon}(x_i)}{\kappa} \leq 1 && \text{if } \bar\eta_i^{\epsilon} = a_i^\epsilon > 0, \\
				%				& - \frac{\bar\varphi^{\epsilon}(x_i)}{\kappa} = 1 && \text{if } \bar\eta_i^{\epsilon} \in (a_i^\epsilon, b_i^\epsilon), 0 <  a_i^\epsilon < b_i^\epsilon,  \\ 
				%				& - \frac{\bar\varphi^{\epsilon}(x_i)}{\kappa} \in [-1,1] && \text{if } \bar\eta_i^{\epsilon} = 0, a_i^\epsilon < 0 < b_i^\epsilon, \\
				%				& - \frac{\bar\varphi^{\epsilon}(x_i)}{\kappa} =- 1 && \text{if } \bar\eta_i^{\epsilon} \in (a_i^\epsilon, b_i^\epsilon),   a_i^\epsilon < b_i^\epsilon < 0, \\
				%				& - \frac{\bar\varphi^{\epsilon}(x_i)}{\kappa} \geq -1 && \text{if } \bar\eta_i^{\epsilon} = b_i^\epsilon < 0, \\
				%				& - \frac{\bar\varphi^{\epsilon}(x_i)}{\kappa} \leq -1 && \text{if } \bar\eta_i^{\epsilon} = a_i^\epsilon < 0, \\
				%			\end{aligned}
			%			\right.
			\quad \text{for } 1 \leq i \leq k \label{eq:OC-1st-variational} 
		\end{align}
	\end{subequations}
	 	Obviously, $\bar y^{\epsilon}$ and $\bar{\bm{\eta}}^{\epsilon}$ satisfy \eqref{eq:OC-1st-state}. Let $\bar{\varphi}^{\epsilon}$ be the unique solution to \eqref{eq:OC-1st-adjoint-state}, which is well-defined as in \Cref{thm:T-classC2}.
	 In order to show \eqref{eq:OC-1st-variational}, we first take any $\bm{\eta} = (\eta_1, \eta_2, \ldots, \eta_k) \in \R^k$ such that $\bm{a} \leq \bm{ \eta} \leq \bm{b}^\epsilon$.
	 Similar to \eqref{eq:1st-OC-direction-active}, there holds
	 \begin{align*}
	 	%\label{eq:1st-OC-direction-active}
	 	0 & \leq T'(\bm{\bar \eta}^{\epsilon})(\bm{\eta} - \bm{\bar \eta}^{\epsilon}) + \kappa [ j(\bm{\eta}) - j(\bm{\bar \eta}^{\epsilon})] 
	 	% \\
	 	= \sum_{i=1}^{k}  \bar{\varphi}^{\epsilon}(x_i)(\eta_i - \bar{\eta}_i^{\epsilon}) + \kappa [ |{\eta}_i| - |{\bar \eta}_i^{\epsilon}| ].
	 \end{align*}
	% for all $\bm{\eta} \in \R^k$ satisfying $\bm{a} \leq \bm{\eta} \leq \bm{b}^\epsilon$.
	A straightforward argument and the second relation in \eqref{eq:ab-epsilon-zero} imply \eqref{eq:OC-1st-variational}.

	\medskip
	\noindent{\it Step 3: Passing to the limit as $\epsilon \to 0^{+}$.}
	The boundedness of $\{\bm{\bar \eta}^\epsilon\}$ in $\R^k$ implies the existence of a subsequence $\{\epsilon_n \}$ such that
	\begin{equation}
		\label{eq:minimizer-control-lim}
		\bm{\bar \eta}^n  \to {\bm{\tilde\eta}} \quad \text{ in } \R^k
	\end{equation}
	with $\bm{\bar \eta}^n := \bm{\bar \eta}^{\epsilon_n}$ for all $n \geq 1$. Consequently, one has from the choice of $\bm{a}$, $\bm{b}^\epsilon$ in \eqref{eq:ab-def}, and the first inequality in \eqref{eq:ab-epsilon-zero} that
	\begin{equation}
		\label{eq:near-minimizer}
		\left\{
		\begin{aligned}
			& \bar\eta_i - \epsilon_0 \leq \tilde{\eta}_i \leq  \bar\eta_i + \epsilon_0 < 4\pi  && \text{if } i \notin I_{4\pi}(\bm{\bar \eta}),\\
			&\bar\eta_i - \epsilon_0 \leq \tilde{\eta}_i \leq 4\pi  && \text{if } i \in I_{4\pi}(\bm{\bar \eta})
		\end{aligned}
		\right.
	\end{equation}
	 for all $1 \leq i \leq k$.
	Setting $\bar y^n := \bar y^{\epsilon_n} = S(\bm{\bar \eta}^n)$ and using the Lipschitz continuity result \eqref{eq:S-continuity-W1q} together with the convergence \eqref{eq:minimizer-control-lim}, we obtain
	\begin{equation}
		\label{eq:state-limit}
		\bar y^n \to \tilde y = S(\bm{\tilde \eta}) \quad \text{strongly in } W^{1,q}_0(\Omega)
	\end{equation}
	for all $1 \leq q < 2$.  
	We now show that $(\tilde y, \bm{\tilde \eta})$ is a global minimizer of \eqref{eq:P-original}. To this end, 
	we first use the continuous embedding $W^{1,q}_0(\Omega) \hookrightarrow L^2(\Omega)$ for $q$ sufficiently close to $2$, and thus deduce from \eqref{eq:state-limit} that
	\begin{equation}
		\label{eq:yn-limit}
		\bar y^n \to \tilde y \quad \text{strongly in } L^2(\Omega).
	\end{equation}
	Combining this with \eqref{eq:minimizer-control-lim} yields
	\begin{equation}
		\label{eq:value-min-limit}
		\hat J(\bm{\bar \eta}^n) \to \hat J(\bm{\tilde \eta}).
	\end{equation}
	Let $(y, \bm{\eta})$ be an arbitrary admissible point  of \eqref{eq:P-original} such that $\bm{ \eta} \neq \bm{\bar \eta}$ and
	\[
		| {\eta}_i - \bar \eta_i| \leq \frac{\epsilon_0}{2}  \quad \text{for all } 1 \leq i \leq k.
	\] 
	By \eqref{eq:SE-existence-solution-assumption-0}, there hold $\bm{\eta}_{\max} \leq 4 \pi$ and $y = S(\bm{\eta})$. 
	Let $\rho \in (0, \frac{\epsilon_0}{2})$ be arbitrary but fixed and let us set $\bm{ \eta}_{\rho} := \bm{ \eta} - \rho \bm{1}$.
	There exists an integer $n_0(\rho) >0$ such that $\bm{a} \leq  \bm{ \eta}_{\rho} \leq \bm{ b }^{\epsilon_n}$ for all $n \geq n_0(\rho)$. 
	Therefore, $\bm{ \eta}_{\rho}$ is a feasible point of  \eqref{eq:P-regularized} corresponding to $\epsilon := \epsilon_n$.
	The optimality in $(\bar y^n,\bm{\bar \eta}^n)$ then gives 
	\[
		\hat J(\bm{\eta} - \rho \bm{1}) \geq \hat J(\bm{\bar \eta}^n) \quad \text{for all } n \geq n_0(\rho).
	\]
	By letting $n \to \infty$ and employing \eqref{eq:value-min-limit}, we have
	$
		\hat J(\bm{\eta}- \rho \bm{1}) \geq \hat J(\bm{\tilde  \eta}).
	$
	Now passing $\rho \to 0^+$ yields $\hat J(\bm{\eta}) \geq \hat J(\bm{\tilde  \eta})$.
	Since $(y, \bm{\eta})$ is arbitrary and $(\bar y, \bm{\bar \eta})$ is a strict local minimizer of \eqref{eq:P-original}, it follows from \eqref{eq:near-minimizer} that $\bm{\tilde \eta} = \bm{\bar \eta}$. 
	
	We have shown that $(\bar y^n, \bm{\bar \eta}^n)$ converges strongly to $(\bar y, \bm{\bar \eta})$ in $W^{1,q}_0(\Omega) \times \R^k$ and $\bar y = S(\bm{\bar \eta})$.
	
	\medskip 
	
	We now show the convergence of $\{\bar\varphi^n := \bar\varphi^{\epsilon_n} \}$ in $H^1_0(\Omega) \cap L^\infty(\Omega)$. To this end, we first deduce from  \eqref{eq:yn-limit} that $\{ \bar y^n\}$ is bounded in $L^2(\Omega)$.
	Moreover, in view of the  estimate in \eqref{eq:exp-S-continuity-L1-positive-part-whole}, there holds
	\begin{equation}
		\label{eq:exp-yn-limit}
		e^{\bar y^n} \to e^{\bar y} \quad \text{strongly in } L^1(\Omega).
	\end{equation}
	Applying now estimate (5) and (6) in \cite[Thm.~1.1]{Lee2025}  to \eqref{eq:OC-1st-adjoint-state} for $\epsilon := \epsilon_n$ yields
	\begin{equation}
		\label{eq:adjoint-state-bounded}
		\norm{\bar{\varphi}^n}_{H^1_0(\Omega)} + \norm{\bar{\varphi}^n}_{L^\infty(\Omega)} \leq C\norm{\bar y^{n} - y_d}_{L^2(\Omega)} \quad \text{for all } n \geq 1
	\end{equation}
	and for some constant $C$ independent of $n$.
	Therefore, $\{\bar{\varphi}^n \}$ is bounded in $H^1_0(\Omega) \cap L^\infty(\Omega)$. 
	There thus exists a subsequence $\{\bar{\varphi}^{n_m} \}$ that converges to some $\tilde{\varphi} \in H^1_0(\Omega) \cap L^\infty(\Omega)$ weakly in $H^1_0(\Omega)$ and weakly-star in $L^\infty(\Omega)$. Letting $m \to \infty$ in the equation \eqref{eq:OC-1st-adjoint-state} for $\bar{\varphi}^{n_m}$ then yields
	\[
		-\Delta \tilde \varphi + e^{\bar y} \tilde\varphi = \bar y - y_d \, \text{in } \Omega, \quad
		\tilde \varphi = 0 \, \text{on } \Gamma,
	\]
%	\[
%	\left\{
%	\begin{aligned}
%		-\Delta \tilde \varphi + e^{\bar y} \tilde\varphi &= \bar y - y_d && \text{in } \Omega, \\
%		\tilde \varphi &= 0 && \text{on } \Gamma,
%	\end{aligned}
%	\right.
%	\]
	where we have used the limits \eqref{eq:exp-yn-limit} and \eqref{eq:yn-limit}. From this and the uniqueness of solutions to \eqref{eq:OC-1st-adjoint-state-P-orig}, one has $\tilde{\varphi} = \bar\varphi$. 
	Consequently, the subsequence $\{\bar{\varphi}^{n_m} \}$  converges weakly in $H^1_0(\Omega)$ and  weakly-star in $L^\infty(\Omega)$ to $\bar\varphi$.
	The uniqueness of $\bar\varphi$ implies that
	the whole sequence 	$\{\bar{\varphi}^n\}$  converges to $\bar\varphi$ weakly in $H^1_0(\Omega)$ and weakly-star in $L^\infty(\Omega)$.
	
	We now show that the convergence  of $\{\bar{\varphi}^n\}$  is indeed in the strong topology of $H^1_0(\Omega)$. For this purpose, subtracting the equations for $\bar{\varphi}^n$ and for $\bar{\varphi}$ yields
	\begin{equation*}
		%\label{eq:adjoint-state-subtraction}
		\left\{
		\begin{aligned}
			-\Delta (\bar{\varphi}^n - \bar \varphi) + e^{\bar y^n} (\bar{\varphi}^n - \bar \varphi) &=\bar y^n - \bar y - (e^{\bar y^n} - e^{\bar y})\bar{\varphi} && \text{in } \Omega, \\
			(\bar{\varphi}^n - \bar \varphi) &= 0 && \text{on } \Gamma.
		\end{aligned}
		\right. 
	\end{equation*}
	Testing this equation by $(\bar{\varphi}^n - \bar \varphi)$ yields
	\[
	\norm{\nabla (\bar{\varphi}^n - \bar \varphi)}_{L^2(\Omega)}^ 2 \leq \int_\Omega [\bar y^n - \bar y - (e^{\bar y^n} - e^{\bar y})\bar{\varphi} ](\bar{\varphi}^n - \bar \varphi) dx.
	\]
	The right-hand side of the last estimate tends to $0$, since $\bar y^n - \bar y - (e^{\bar y^n} - e^{\bar y})\bar{\varphi} \to 0$ strongly in $L^1(\Omega)$ (see \eqref{eq:yn-limit} and \eqref{eq:exp-yn-limit}) and $\{\bar{\varphi}^n\}$ weakly-star converges to $\bar \varphi$ in $L^\infty(\Omega)$. 
	Therefore, $\{\bar{\varphi}^n\}$  converges to $\bar\varphi$ strongly in $H^1_0(\Omega)$.
	
	%It remains to prove \eqref{eq:OC-1st-variational-P-orig}.
	The final step is to prove \eqref{eq:OC-1st-variational-P-orig} and the continuity around $x_i$ of $\bar\varphi$. Assume now that $\bar \eta_i < 4\pi$. We shall show that $\bar\varphi$ is continuous in a neighborhood of $x_i$. Indeed, by assertion \ref{item:less-4pi} in \Cref{prop:critical-4pi}, there exist constants $\rho >0$ and $\tau >0$ such that
	\[
	e^{\bar y} \mid_{B_{\R^2}(x_i, \rho)} \in L^{1 + \tau}(B_{\R^2}(x_i, \rho)).
	\]
	On the other hand, it is easy to see that 
	\[
	-\Delta \bar\varphi = - e^{\bar y} \bar\varphi + \bar y - y_d \quad \text{in } \mathcal{D}'(B_{\R^2}(x_i, \rho)).
	\]
	From this and the interior H\"{o}lder regularity of solutions to Poisson's equation; see, e.g., \cite[Thm.~8.24]{Gilbarg_Trudinger}, we have the continuity around $x_i$ of
	$\bar\varphi$. This, together with the limit $\bar{\varphi}^n \to \bar{\varphi}$ in $H^1_0(\Omega)$, implies that
	\[
		\bar{\varphi}^n(x) \to \bar{\varphi}(x) \quad \text{for all } x \in  B_{\R^2}(x_i, \rho_1)
	\]
	for some $\rho_1 >0$. 
	From this, together with \eqref{eq:OC-1st-variational} and the definitions of  $\bm{a}$ and $\bm{b}^\epsilon$ in \eqref{eq:ab-def}, we obtain \eqref{eq:OC-1st-variational-P-orig} by applying a standard argument.
\qed

%%%% Second-order OCs when $\omega_{\max}$ < 4\pi. It is a direct consequence of \Cref{thm:2nd-OCs-nec}

%\todo{define the critical cone and $T''$?}
%\begin{theorem}
%	\label{thm:2nd-OCs-nec-P-orig}
%	Assume that $\bar{\bm{\eta}}$ is a local minimizer of \eqref{eq:P-original} such that 
%	\begin{equation}
	%		\label{eq:minimizer-less-4pi-P-orig}
	%		\bar{\bm{\eta}}_{\max} < 4 \pi.
	%	\end{equation}
%	Then, there holds
%	\begin{equation}
	%		\label{eq:2nd-OCs-nec-P-orig}
	%		T''(\bm{\bar \eta})\bm{\eta}^2 \geq 0 \quad \text{for all } \bm{\eta} \in C_{\bar{ \bm{\eta}}}.
	%	\end{equation}
%\end{theorem}

\section{Conclusions}
We have investigated an optimal control problem governed by a semilinear elliptic equation with exponential nonlinearity, where the source term includes a linear combination of Dirac measures concentrated at a finite set of distinct points.
Since the state equation admits a unique solution only when the masses of the Dirac measures do not exceed the threshold value $4\pi$, the control-to-state operator is continuously differentiable only on an open subset of the control space.
To derive first-order optimality conditions, we considered a family of regularized optimal control problems. These regularizations incorporated suitable control constraints, ensuring that the admissible controls remained within the domain of differentiability of the control-to-state operator. By passing to the limit in these regularized problems, we obtained an optimality system for the original problem.

\begin{acknowledgements}
	This work was partly completed during a visit of the	author to Vietnam Institute for Advanced Study in Mathematics (VIASM). The author would
	like to thank VIASM for their financial support and hospitality.
\end{acknowledgements}

%\appendix  %This command ends the counting of sections.
%\section*{Appendix:  Instructions for Appendices}

\appendix
\section{Estimates of Green's function on the unit ball}

\begin{lemma}
	\label{lem:Green-esti-unitball}
	Let $G$ be the Green function of the Laplace operator $-\Delta$ on the unit ball in $\R^2$. Then, there holds
	\begin{equation}
		\label{eq:Green-esti-unitball}
		0 \leq G(x,y) \leq \frac{1}{2 \pi} \ln \left( \frac{2}{|x-x'|}\right) \quad \text{for all } x, x' \in {B_{\R^2}(0,1)}, \quad x \neq x'. 
	\end{equation}
\end{lemma}
{\it Proof.}  
	We have for all $x, x' \in B_{\R^2}(0,1)$ with $x \neq x'$ that
	\[
	G(x,x') = 
	\left\{
	\begin{aligned}
		& - \frac{1}{2  \pi} \ln |x'| && \text{if } x = 0,\\
		& - \frac{1}{2 \pi}[ \ln {|x-x'|} - \ln ( | \frac{x}{|x|} - |x| x' | )   ] && \text{if } x \neq 0;
	\end{aligned}
	\right.
	\]
	see, e.g. \cite[Form.~(2.23), Chap.~2]{Gilbarg_Trudinger} and \cite[Form.~(41) \& (6), Chap.~2]{Evans2010}.
	The nonnegativity of $G$ follows from \cite[Form.~(2.24), Chap.~2]{Gilbarg_Trudinger} with noting that the Green function defined in \cite{Gilbarg_Trudinger} is associated with the $\Delta$ operator. 
	However, we now state the detailed argument showing both two estimates in \eqref{eq:Green-esti-unitball}. To this end, we first consider the case $x = 0$ and have
	\[
	0 \leq G(0, x') =  - \frac{1}{2  \pi} \ln |x'| =   \frac{1}{2  \pi} \ln \frac{2}{|0-x'|},
	\]
	which shows \eqref{eq:Green-esti-unitball}. 
	For $x \neq 0$, we have
	\begin{align*}
		| \frac{x}{|x|} - |x| x' |^2 & = 1 - 2 x \cdot x' + |x|^2|x'|^2  = |x-x'|^2 + 1 - |x|^2 - |x'|^2 + |x|^2|x'|^2\\
		& = |x- x'|^2 + (1-|x'|^2)(1- |x|^2),
	\end{align*}
	which thus gives
	$
	|x-x'| \leq | \frac{x}{|x|} - |x| x' | \leq 2
	$
	for all $x, x' \in B_{\R^2}(0,1)$. \eqref{eq:Green-esti-unitball} then follows.
\qed   

%%%% the openess of
%\section{Openness of subsets \texorpdfstring{{\boldmath $\mathcal{D}_{ \bm{b} }$}}{}  and \texorpdfstring{$V_{\bm{b} } $}{}} % defined in  \eqref{eq:Db-set} and \eqref{eq:Vb-set} }% $\mathcal{D}_{ \bm{b} }$ and $V_{\bm{b} } $ }
\section{Openness of subsets {{\boldmath $\mathcal{D}_{ \bm{b} }$}}{}  and {$V_{\bm{b} } $}{}} % defined in  \eqref{eq:Db-set} and \eqref{eq:Vb-set} }% $\mathcal{D}_{ \bm{b} }$ and $V_{\bm{b} } $ }
\label{sec:openness-Db-Wrb}

\begin{lemma}
\label{lem:openness-Db-wWrb}
Assume that $\bm{b} \in \R^k$ and $p>1$ such that $r_{\bm{b}}$ determined in \eqref{eq:rb-constant} is greater than $1$. 
Let  $\mathcal{D}_{\bm{b}}$, $V_{\bm{b}}$, $\mathcal{D}$, and $W_{r_{\bm{b}}}$ be defined as in  \eqref{eq:Db-set}, \eqref{eq:Vb-set}, \eqref{eq:D-space}, and \eqref{eq:Wrb-space}, respectively. Then, $\mathcal{D}_{\bm{b}}$ is open in $\mathcal{D}$ and $V_{\bm{b}}$ is open in $W_{r_{\bm{b}}}$.
\end{lemma}
{\it Proof.}  
We first prove the openness of $\mathcal{D}_{\bm{b}}$ in $\mathcal{D}$. To this end, for any $\mu_* := \sum_{i =1}^k \omega_i^* \delta_{x_i} \in \mathcal{D}_{\bm{b}}$ with ${\omega}^*_i < b_i$ for all $1 \leq i \leq k$, choosing $\epsilon >0$ such that 
\[
\epsilon < \frac{1}{2}\min \{ b_i - \omega_i^* \mid 1 \leq i \leq k \}
\]
yields that the open ball in $\mathcal{D}$ centered at $\mu_*$ with radius $\epsilon$ belongs to $\mathcal{D}_{\bm{b}}$.

We now show the openness of $V_{\bm{ b }}$ in $W_{r_{\bm{b}}}$. For that purpose, let $y_* \in V_{\bm{ b }}$ be arbitrary, but fixed. By definition, there exists a pair $(\mu_*, h_*) \in \mathcal{D}_{\bm{b}} \times L^{r_{\bm{b}}}(\Omega)$ such that
\[
- \Delta y = \mu_* + h_* \quad \text{in } \Omega, \quad y_* \in W^{1,1}_0(\Omega).
\]
The openness of $\mathcal{D}_{\bm{b}}$ in $\mathcal{D}$ guarantees the existence of $\epsilon >0$ such that
$B_{\mathcal{D}}(\mu_*,\epsilon) \subset \mathcal{D}_{\bm{b}}$. Taking any $y \in W_{r_{\bm{b}}}$ such that $\norm{y - y_*}_{W_{r_{\bm{b}}}} < \epsilon$ gives the existence of a unique pair $(\mu, h) \in \mathcal{D} \times L^{r_{\bm{b}}}(\Omega)$ satisfying
\[
- \Delta y = \mu + h \quad \text{in } \Omega, \quad y \in W^{1,1}_0(\Omega).
\]
Moreover, there holds
\[
\norm{\mu - \mu_*}_{\mathcal{M}(\Omega)}  \leq \norm{y - y_*}_{W_{r_{\bm{b}}}} < \epsilon.
\]
This indicates that $\mu \in B_{\mathcal{D}}(\mu_*,\epsilon) \subset \mathcal{D}_{\bm{b}}$ and thus $y \in V_{\bm{b}}$. 
%The proof of the lemma is thus complete.
\qed   

\bibliographystyle{spmpsci}
\bibliography{exponentialnonlinearity_fulldois}

\end{document}